\documentclass[aap,MSNbibl,seceqn,dvips]{arximspdf}
\usepackage{mathbh}

\doi{10.1214/12-AAP869} 
\volume{23}
\issue{3}
\pubyear{2013}
\firstpage{1219}
\lastpage{1253}

\makeatletter

\def\varnothing{\varnothing}
\def\ti{\to\infty}

\def\e{\varepsilon}

\newcommand{\R}{\mathbb{R}}

\def\1{\mathbh{1}}
\newcommand{\ssup}[1]{({#1})}

\newtheorem{theorem}{Theorem}[section]
\newtheorem{Lemma}[theorem]{Lemma}
\newtheorem{cor}[theorem] {Corollary}
\newproclaim{remark}[theorem] {Remark}

\newcommand{\eps}{\varepsilon}

\newcommand{\Fcal}{\mathcal{F}}
\newcommand{\Gcal}{\mathcal{G}}

\newcommand{\Xcal}{\mathcal{X}}
\renewcommand{\e}{\operatorname{e}}
\newcommand{\df}{{\,\stackrel{\mathrm{def}}{=}\,}}
\newcommand{\dff}{{\,\stackrel{def}{=}\,}}

\newcommand{\eqref}[1]{(\ref{#1})}

\makeatother

\begin{document}
\begin{frontmatter}

\title{On a preferential attachment and generalized P\'olya's urn model}
\runtitle{\hspace*{-10pt}Preferential attachment and generalized P\'olya's urn model}

\begin{aug}
\author[A]{\fnms{Andrea} \snm{Collevecchio}\corref{}\ead[label=e1]{collevec@unive.it}\thanksref{t1}},
\author[B]{\fnms{Codina} \snm{Cotar}\ead[label=e2]{cotar@ma.tum.de}\thanksref{t1}}
\and
\author[C]{\fnms{Marco} \snm{LiCalzi}\ead[label=e3]{licalzi@unive.it}}
\thankstext{t1}{Supported by IAS-Technische Universitaet Muenchen.}
\runauthor{A. Collevecchio, C. Cotar and M. LiCalzi}
\affiliation{Universit{\`a} Ca' Foscari Venezia and Monash University,
Technische~Universitaet Muenchen, and Universit{\`a} Ca' Foscari Venezia}
\address[A]{A. Collevecchio\\
Department of Management\\
Universit{\`a} Ca' Foscari Venezia\\
Italy\\
and\\
Department of Mathematics\\
Monash University\\
Melbourne\\
Australia\\
\printead{e1}} 
\address[B]{C. Cotar\\
Zentrum Mathematik\\
Technische Universitaet Muenchen\\
Germany\\
\printead{e2}}
\address[C]{M. LiCalzi\\
Department of Management\\
Universit{\`a} Ca' Foscari Venezia\\
Italy\\
\printead{e3}}
\end{aug}

\received{\smonth{8} \syear{2011}}
\revised{\smonth{4} \syear{2012}}

%
\begin{abstract}
We study a general preferential attachment and P\'olya's urn model. At
each step a new vertex is introduced, which can be connected to at most
one existing vertex. If it is disconnected, it becomes a pioneer
vertex. Given that it is not disconnected, it joins an existing pioneer
vertex with probability proportional to a function of the degree of
that vertex. This function is allowed to be vertex-dependent, and is
called the reinforcement function. We prove that there can be at most
three phases in this model, depending on the behavior of the
reinforcement function. Consider the set whose elements are the
vertices with cardinality tending a.s. to infinity. We prove that this
set either is empty, or it has exactly one element, or it contains all
the pioneer vertices. Moreover, we describe the phase transition in the
case where the reinforcement function is the same for all vertices. Our
results are general, and in particular we are not assuming monotonicity
of the reinforcement function.

Finally, consider the regime where exactly one vertex has a degree
diverging to infinity. We give a lower bound for the probability that a
given vertex ends up being the leading one, that is, its degree
diverges to infinity. Our proofs rely on a generalization of the Rubin
construction given for edge-reinforced random walks, and on a Brownian
motion embedding.
\end{abstract}

%

\end{frontmatter}

\section{Introduction}\label{intro}

\subsection{Setting and {m}otivation}

We study the following model. Given finitely many classes (or groups)
each containing a given initial number of members, new members arrive
one at a time. For each new member arriving at time $n$, with
probability\vadjust{\goodbreak} $s_n\ge0$ we create a new class in which we place the
member; with probability $1-s_n$, we place the member in an existing
class. We assume that each existing class attracts new members with
probability proportional to a certain positive function of the
cardinality of the group, called the \textit{reinforcement} or \textit
{weight scheme} $f$. If the groups are allowed to have different
reinforcement schemes, then we show that looking at the asymptotics as
time tends to infinity we have exactly three different regimes: one
group is infinite and all the others are finite; all groups are
infinite; all groups are finite. Our main result,
Theorem~\ref{princres1}, shows that in the first regime the process
will eventually create a unique infinite group: this happens when each
group is reinforced quite a bit, but not too much with respect to the
other groups. In the second regime, the cardinality of each group goes
to infinity. Finally, in the last regime, all the groups will be
finite; what happens is that the process creates various peaks: in the
beginning one group dominates the others, but sooner or later another
group will start dominating, and this change happens infinitely many
times. In this way, no group dominates definitively the other groups.
This is a kind of ``there is always a faster gun'' principle.

Our model is a generalization of two models from two different classes:
one model from the class of \textit{preferential attachment models}, as
introduced in~\cite{DR2} and in~\cite{KS}, and one model from the class
of \textit{reinforcement processes}, as introduced in~\cite{FC}.

The first main model we are generalizing was introduced and studied
independently in~\cite{DR2} and in~\cite{KS}, and later studied in more
detail in~\cite{ROJS} and~\cite{Ru}. This model is part of the class of
\textit{preferential attachment models}, which are models of growing
networks, and which were first proposed in the highly-influential
papers~\cite{BA1} and~\cite{BA2}. In~\cite{BA1} new vertices arrive at
the network one at a time and send a
fixed number $m$ of edges to already existing vertices; the probability
that a new vertex is linked to a given existing vertex is
proportional to the in-degree of the respective existing vertex. Here,
the in-degree of a vertex is the number of children of that vertex.

The model studied in~\cite{DR2,KS,ROJS} and~\cite{Ru}
is as follows: consider a model of an evolving network in which
new vertices arrive one at a time, each connecting by an edge to a
previously existing vertex with a probability proportional to a
function~$f$ of the existing vertex's in-degree. This function $f$ is
called \textit{attachment rule}, or \textit{weight function}, and it
determines the existence of two main different regimes. The first
regime corresponds to $f(j)=j+1$, and it was studied in~\cite{BA1,BA2} and~\cite{Ru}; the second regime corresponds to for $
\gamma<1$, and it was studied in~\cite{Ru}. The third regime
corresponds to $f(j)=(j+1)^\gamma$ for $ \gamma>1$, and it was
studied in~\cite{ROJS}. In the first two regimes, it is shown that the
degrees of all vertices grow to infinity; in the third regime there is
a second phase as one vertex eventually dominates all other vertices.
In the first regime, the so-called P\'olya urn, the urn process is
exchangeable and is the only case where exchangeability appears; see
\cite{HLS}.
(For more results on preferential attachment models, see the survey
\cite{Bh}.)\vadjust{\goodbreak}

Preferential attachment models have been motivated by real-life
problems, especially in regards to network and internet applications.
One important example of growing networks is the World Wide Web, in
which the more popular a page (or vertex) is the more hits it
{receives}; a similar principle applies to social interaction or to
citation networks. Another example is the one of users of a software
program who can report bugs on a website. Bugs with the highest number
of requests get priority to be fixed. If the user cannot find
an existing report of the bug, they can create a new report. However,
it could be that there are duplicate reports, in which case the number
of requests
is split between the reports, making it less likely that the bug the
user found will get fixed. Since bugs that have more
requests appear higher up the search results, the user is more likely
to add a request to an existing report than to a new one.

This can be explained by the fact that such networks are built
dynamically and that new vertices prefer to attach themselves to
existing popular vertices with high in-degree rather than to existing
unpopular vertices with low in-degree.

The second main model we are generalizing is studied in~\cite{FC,oliv} and~\cite{tzhou09}. It is known as the \textit{generalized
P\'olya urn process}; it belongs to the class of \textit{reinforcement
processes} and can be described as follows.
Given finitely many bins each containing one ball, new balls arrive
one at a time. For each new ball, with probability $p\ge0$ we create
a new bin in which we place the ball; with probability $1-p$, we place
the ball in an existing bin. The probability that the ball is placed in
an already existing bin is proportional to $f(j)=j^\gamma$, where $j$ is
the number of balls in that bin. The case with $p=0$ and $\gamma=1$ is
the well-known \textit{P\'olya urn problem}. For $p=0$ and $\gamma>0$
no new bins are
created, and the process is called a \textit{finite P\'olya process
with exponent $\gamma$}. If $p>0$, then the process is called an \textit
{infinite P\'olya process}. Similar{ly} to the preferential attachment
models, for generalized P\'olya urn processes {with} $f(j)=j^\gamma$,
it is {known that for} $\gamma\le1$ the number of balls in all bins
eventually grows to infinity, whereas for $\gamma>1$ one bin eventually
{comes to} dominate all other bins. (A detailed review of a number of
other interesting results on P\'olya's urn processes and on
reinforcement processes in general is provided in the survey~\cite{Pem}.)

The generalized P\'olya urn process has applications to many areas. We
briefly mention one such application to biology; for an extensive
overview of other applications of generalized P\'olya urn processes to
reinforced random walks, statistics, computer science, clinical trials,
biology, psychology and economics, see, for example, Chapter 4 in~\cite{tzhou09}.

The generalized P\'olya urn process with $p=0$ is used in~\cite{DAGP}
and~\cite{SKJC} to study a real-life application; the reinforcement
scheme use{d} in these papers is set to $f(j) = j^\gamma${, with
}$\gamma>0$, and
real-life data {are} compared {against} different values of $\gamma$
and initial configurations.
More precisely, the authors study a colony of ants, which explores a
chemically unmarked territory randomly, starting from the nest. The exploration
is done on a fixed number $k$ of paths of various lengths. Each ant
passes along one of the $k$ paths leaves a pheromone mark and in
this way infuences the
following ant's decision in choosing a particular path. This decision
is also influenced by whether the paths of various {lengths} are
discovered at the same time, or whether they are discovered at
different times. In the real-life experiment it is noticed in the case
of paths of equal lengths that, after initial fluctuations, one of the
paths becomes more
or less completely preferred to the others.

We will show in our paper that the above two models, belonging to
these two different areas, are in fact closely related because they are
both special cases of our much more general model. The first of our
results, Theorem~\ref{princres}, proved for our general model, unifies
the two above-described phase transition results for a very general
class of weight functions $f$; the result holds in particular both for
preferential attachment processes and for generalized P\'olya's urn
processes. It is worth noting that our condition on the weight function
is much weaker than all previously-proved results for the models we
generalize. Moreover, in our main result, Theorem~\ref{princres1}, we
show, under no assumptions on the weight function, that we can have
only three possible phases; in the third phase, all groups (resp.,
vertices, bins) stay finite as time tends to infinity. To the best of
our knowledge, this is the first time when a third regime as described
in our Theorem~\ref{princres1}, has been proved for any model of
preferential attachment or P\'olya's urn type. In the case of weight
functions $f$ which give rise to the second phase, we devise in our
Theorem~\ref{testo1}, and, respectively, in Corollary~\ref{rf}, a test
for obtaining an upper bound, and, respectively, a lower bound, on the
probability that a given group ends up being {dominant}.

The motivation for our model comes from the class of \textit{species
sampling sequences}, to which class our model belongs. Species sampling
sequences are models for exchangeable sequences $(X_n)$ with a
prediction rule, that is, a formula for the conditional distribution of
${X_{n+1}}$ given $X_1,X_1,\ldots, X_n$ for $n=1,2,\ldots,n.$ More
precisely, given the first $n$ terms of the sequence $(X_n)$, $X_{n+1}$
equals the $i$th distinct value observed so far with probability
$p_{n+1,i}$, for $i=1,2,\ldots, n$, and otherwise $X_{n+1}$ is a new
value with distribution $\nu$ for some probability measure $\nu$.
Species sampling sequences were first introduced and studied in \cite
{pit1,pit2} and are now used extensively in Bayesian
nonparametric statistics; see, for example,~\cite{hp,jam} or
\cite{lqmt} for more on species sampling sequences or for their
applications to statistics.

We next introduce precisely our model.

\subsection{The main model}\label{MAIN}
We consider the following model where at each step a new vertex and at
most one new edge appear according to the following rules. The
probability that the new vertex is disconnected is positive and may
change in time. When a vertex is disconnected from the existing ones,
it becomes a \textit{pioneer vertex}. We label the pioneer vertices in
order of appearance. Given that the new vertex is connected to an
existing one, the latter is chosen with probability proportional to a
reinforcement scheme of its degree. The graph formed with this
procedure is the union of trees. Each tree has a pioneer vertex as a
root. The tree with root $i$ observed at time $n$, is called the $i$th
group (or $i$th component) by time $n$.

More formally, fix a collection of positive functions $f_{k}
\dvtx\mathbb N \to\R^{+}$ with $f_{k}(0) = 0$ and $f_{k}(i)>0$ for
all $i,k \ge1$, {and a sequence $\{s_{n}\}$ which takes values in
$[0,1]$}. Set $A_{1}(1)=1$ and $A_{j}(1) = 0$ for all $j \ge2$. Set
$L_{1} =1$. We define the random variables $A_{i}(n+1)$ and $L_{n+1}$
recursively as follows:
\begin{eqnarray}
\mathbb{P} \bigl(A_{k}(n+1) = A_{k}(n) +1 | L_{n},
\bigl\{A_{j}(n), \mbox{ with $j \in\mathbb{N}$}\bigr\} \bigr) &=& (1-
s_n )\frac
{f_{k}(A_{k}(n))}{\sum_{s=1}^{L_{n}} f_{s}(A_{s}(n))}\nonumber\\
\eqntext{\mbox{for $ i \le L_{n}$},}
\\
\mathbb{P}\bigl(A_{L_{n}+1}(n+1) = 1 | L_{n}, \bigl
\{A_{j}(n), \mbox{with $j \in\mathbb{N}$}\bigr\} \bigr) &=&
s_{n},\nonumber
\end{eqnarray}
while $A_{j}(n+1) = 0$ for all $j >L_{n}+1$. Moreover,
\[
L_{n+1} \df\max\bigl\{ j \ge1 \dvtx A_{j}(n+1) \ge1\bigr\}.
\]

Notice that $A_{i}(n+1) - A_{i}(n) \in\{0,1\}$, and exactly for one
index $i$ this difference equals 1. The random variable $A_{i}(n)$ is
the cardinality of the $i$th group by time~$n$.
We call the process $\{A_{i}(n), i,n \ge1\}$ a \textit{generalized
attachment model} whose parameters are the sequence $\{s_{n}\}$ and the
reinforcement functions $\{f_{k}\}$, abbreviated with $\operatorname{GAM}(\{f_{k}\},\{
s_{n}\})$. We emphasize the fact that we do not make any assumptions on
the update functions $\{f_{k}\}$, other than positiveness, and $\{
s_{n}\}$ is allowed to be random. As shown in Theorem~\ref{princres1},
some of our strongest results hold for a group-dependent deterministic
reinforcement scheme $\{f_{j}\}$, that is, where each group $j$ follows
its own reinforcement scheme $f_j$, independently of the other groups.
From the point of view of applications, this allows one to take into
account the case where different groups have different update schemes,
which is what would be expected in many real-life situations. We use
the symbol $\operatorname{GAM}(f, \{s_{n}\})$ to denote a
generalized attachment model where the update functions $f_{k}$ are
equal to the positive function $f$ for each $k \ge1$.


We briefly discuss next the link of our work to the recent literature.
The two main models that we generalize were studied in detail in the
particular case with reinforcement scheme proportional to $f(j)=j^\gamma
$, where $\gamma>0$.


Let us look first at the literature on preferential attachment models
connected to our generalized attachment model. The preferential
attachment model studied in~\cite{DR2,KS,ROJS} and \cite
{Ru} is just $\operatorname{GAM}(f,\{s_{n}\})$ for the particular case of
$s_n=m(n)c/ (\sum_{s=1}^{L_{n}} f(A_{s}(n))+m(n) c )$, where we
denoted by $m(n)\le n$ the number of groups (resp., vertices) with no
children at time $n$, and where $c>0$. Then in the growing network,
$A_{j}(n)\ge1$ represents the in-degree at time $n$ of existing vertex
$j$ with strictly positive in-degree; that is, vertex $j$ has
$A_{j}(n)$ children. With probability $f(A_{j}(n)) / (\sum_{s=1}^{L_{n}} f(A_{s}(n)+m(n) c  )$, a new arriving vertex attaches
to an existing vertex $j$ with strictly positive in-degree $A_{j}(n)$;
with probability $m(n)c/ (\sum_{s=1}^{L_{n}} f(A_{s}(n))+m(n) c
)$, a new arriving vertex attaches to one of the existing $m(n)$
vertices with $0$ in-degree. For the case where the reinforcement
function $f$ is linear or super-linear, $\{s_n\}$ is bounded away from
$1$, so we can apply our results to the case of preferential attachment models.

In~\cite{ROJS} the authors look at the preferential attachment model
with reinforcement scheme $f(j)=(j+1)^\gamma,\gamma>1,$ for which they
prove a similar result to our Theorem~\ref{princres} by using the
original Rubin construction. In~\cite{BOR} and~\cite{mori},
respectively, in~\cite{Ru}, the authors give
the {limiting} degree distribution for a wide range of linear,
respectively, sub-linear, weight functions.

A different preferential attachment model was studied in~\cite{peter1}
and~\cite{peter2}. In this model a new vertex arrives at each step and
attaches to every existing vertex independently with a
probability proportional to a concave weight function $f$ of the
existing vertex's degree. In~\cite{peter1} the authors prove in Theorem
1.5 the same type of phase transition as in our Theorem~\ref{princres},
and they study the degree distribution. In~\cite{peter2} they study the
existence of a giant component, that is, of a connected component
containing a positive fraction of all vertices.

We turn now to the literature on the generalized P\'olya urn model.
This model corresponds to $\operatorname{GAM}(f,\{s_{n}\})$ in the particular case
with $s_n\equiv p$ for all $n\ge n_0$, for some fixed $n_0\in\mathbb
{N}$. In~\cite{FC} the authors consider both the generalized P\'olya
urn model with $p=0$, when the number of bins is fixed, and with $p>0$,
and they prove by combinatorics techniques a similar result to our
Theorem~\ref{princres} for the case of power functions. The case with
$p>0$ and $\gamma<1$ is studied in~\cite{FC} under two additional
assumptions involving the power function $f(j)=j^\gamma$, assumptions
whose validity is left as an open problem in that paper.

The generalized P\'olya urn model with $p=0$ was also the main object
of study in~\cite{oliv} and~\cite{tzhou09}. In~\cite{oliv} the author
studied the case of two fixed bins under a number of technical
assumptions on the function $f$, which exclude, for example, the
(super)-exponential functions, and which assumptions are stated in
Section 4 of that paper.
Theorem 3.3.1 {in}~\cite{tzhou09} proves a result similar {to our}
Theorem~\ref{princres}(i) for the case of a fixed number $m$ of bins
and under the assumption of monotonicity on the super-linear function $f$.

Last, we provide below a definition of species sampling sequences and
why $\operatorname{GAM}(f,\{s_{n}\})$ is such a sequence. Consider a Polish space
$\Xcal$, and let $\mu(\cdot)$ be a diffuse probability measure on $\Xcal
$, that is, $\mu(\{x\}) =0$, for all $x \in\Xcal$. Denote with $\1_{A}$ the indicator function\vadjust{\goodbreak} of the event $A$. A sequence of random
variables $X_{n}$, with $n \ge1$, on $\Xcal$ which has the distribution
%
\begin{equation}
\label{uno} \mathbb{P} (X_{n+1} \in B | X_{1}, \ldots,
X_{n}) =\sum_{i=1}^{n}
p_{n+1,i} \1_{\{ X_{i} \in B\}} + r_{n+1} \mu(B),
\end{equation}
is called a \textit{species sampling sequence} whenever $ r_{n} + \sum_{i}
p_{n,i} =1$, $r_{n}, p_{n,i} \ge0$, and $r_{n}, p_{n,i}$ are $\Fcal_{n-1}$ measurable, where $\Fcal_{n} = \sigma\{X_{1}, X_{2}, \ldots,
X_{n}\}$. It corresponds to $\operatorname{GAM}(f,\{s_{n}\})$ for the case with
$p_{n,i}=f_{i}(A_{i}(n))/\sum_{s=1}^{L_{n}} f_{s}(A_{s}(n))$ and
$r_n=s_n$ for all $n\ge1$. In particular, the Blackwell--MacQueen urn
scheme, also known as Chinese restaurant process, is a species sampling
sequence with the choice $s_{n+1} = p_{n+1,i}= 1/(1+n)$; it
corresponds to $\operatorname{GAM}(f,\{s_{n}\})$ with $f(j)=1/j$ and $s_n=r_n=1/n$ for
all $n\ge1$.

In this paper we give a complete characterization of the existing
phases for a very general class of update functions, for the case
$f_{j} \equiv f$ and for all nonnegative random sequences $(s_n)_{n\in
\mathbb N}$, with $s_n\le p<1$ for all $n\ge1$. In particular we do
not assume any monotonicity on $f$, and our only assumption on $f$ is
for Theorem~\ref{princres}(i), and it controls the oscillation of the
reinforcement function. Moreover, we prove in Theorem~\ref{princres1}
that for any group-dependent
deterministic reinforcement scheme $\{f_{j}\}$, where $\{f_{j}\}$ are
only assumed to be positive, we can only have three possible phases. We
prove the existence of a third phase by an example.
We emphasize the fact that exactly three phases are admitted for this
model.\vspace*{-2pt}

\subsection{Results}

The following are our main results.\vspace*{-2pt}

\begin{theorem}\label{princres} Consider a $\operatorname{GAM}(f,\{s_{n}\})$. Suppose
that $s_{n} \le p$, for some $p<1$ and all $ n \ge1$. Let
\[
\theta^2_{k} =1\Big/ \Biggl(\sum_{s=k+1}^{\infty}
\frac{2}{(1-p)^2f^2(s)} \Biggr).\vspace*{-2pt}
\]
\begin{longlist}[(ii)]
\item[(i)] If
%
\begin{equation}
\label{spb} %
\begin{aligned} \sum_{n=1}^{\infty}
\frac1{f(n)} < \infty\quad \mbox{and}\quad  \sum_{k=1}^{\infty}
\prod_{s=1}^{k} \frac{1}{1 +  (\theta_{k}/(f(s)(1-p)) )}< \infty,
\end{aligned} %
\end{equation}
then there will be, a.s., exactly one group whose cardinality tends to
infinity, all the other groups being finite.
\item[(ii)] The cardinality of each (created) group tends to infinity
a.s. if and only if
\[
\sum_{n=1}^{\infty} \frac1{f(n)} = \infty.
\]
\end{longlist}
\end{theorem}
\begin{remark}
{The second condition in (\ref{spb}) describes a large class of
sequences $f(i)$ whose reciprocal is summable.\vadjust{\goodbreak} In particular it
contains all the} {monotone sequences}, and all the convergent series
with $f(k)\ge k, k\ge1$. This condition is used to control the
oscillation in the sequence. We believe that only the first condition
in (\ref{spb}) is needed for the conclusion.

We will show that all the monotone sequences satisfy the second
assumption in~(\ref{spb}); to prove (\ref{spb}) for $f(k)\ge k$ follows
similar reasoning and will be omitted. We have for large $k\ge k_0>0$,
\begin{eqnarray*}
\sum_{s=k+1}^{\infty} 1/f^2(s)&=&\sum
_{s=k+1}^{k+1+[f(k+1)]} 1/f^2(s)+\sum
_{s=k+1+[f(k+1)]}^{\infty} 1/f^2(s)\le\frac
{f(k+1)}{f^2(k+1)}
\\
&&{}+\frac{1}{f(k+1+[f(k+1)])}\sum_{s=k+1+[f(k+1)]}^{\infty} 1/f(s)
\\
&\le&\frac{f(k+1)}{f^2(k+1)}+\varepsilon\frac{1}{f(k+1)}=\frac
{1}{f(k+1)}
\frac{\varepsilon+1}{\varepsilon}
\end{eqnarray*}
for some $\varepsilon>0$, where for the last inequality we used the fact
that\break $\sum_{s=k+1+[f(k+1)]}^{\infty} 1/f(s)$ converges to $0$ as
$n\rightarrow\infty$. We denoted by $[f(k+1)]$ the integer part of $f(k+1)$.
Therefore
\begin{eqnarray*}
\sum_{k=1}^{\infty} \prod
_{s=1}^{k} \frac{1}{1 +  (\theta_{k}/(f(s)(1-p)) )} &\le& C_{1}+
\sum_{k=k_0}^{\infty} \prod
_{s=1}^{2} \frac{1}{1 +  (\theta_{k}/(f(s)(1-p)) )}\\
&\le& C_{1}+
C_{2}\sum_{k=k_0}^{\infty}
\frac{1}{f(k+1)}.
\end{eqnarray*}
\end{remark}


\begin{remark}
If we remove the hypothesis that $s_{n}$ is bounded away from one, and
suppose that $\sum_{n=1}^{\infty} (1-s_{n})<\infty$, then
by Borel--Cantelli's lemma there exists a random time $N$ such that
for any time $n \ge N$ a new group is formed. Hence the cardinality of
each group will remain finite, and only finitely many groups will end
up having a cardinality larger than $1$. We do not study the case of
$\limsup_{n \ti} s_{n}=1$ and $\sum_{n=1}^{\infty} (1-s_{n})=\infty$.
\end{remark}
The following result is a corollary to the proof of Theorem \ref
{princres}(i). It generalizes the results contained in~\cite{ROJS}
about the degree of vertices in the preferential attachment model.
%
\begin{cor}\label{findegree} If the hypothesis of Theorem \ref
{princres}\textup{(i)} holds, then\break $\lim_{n \ti} A_{i}(n)>1$ for only finitely
many $i$.
\end{cor}
The following theorem establishes that $\operatorname{GAM}(\{f_{j}\},\{s_{n}\})$ can
have only three possible phases. The theorem holds true if the $f_{j}$
are random functions independent of $s_{n}$ satisfying the conditions
of the theorem almost surely.\vadjust{\goodbreak}

\begin{theorem}\label{princres1} Consider a $\operatorname{GAM}(\{f_{j}\},\{s_{n}\})$.
Suppose that $ s_{n} \le p <1$, for some $p<1$ and all $ n \ge1$.
\begin{longlist}[(ii)]
\item[(i)] If
\[
\sum_{n=1}^{\infty} \frac1{f_{j}(n)}
< \infty \qquad\mbox{for at least one {created group} $j \in\mathbb{N}$},
\]
then there will be, a.s., \textup{at most} one group whose cardinality
tends to infinity, all the other groups being finite.
\item[(ii)] If
\[
\sum_{n=1}^{\infty} \frac1{f_{j}(n)}
= \infty\qquad \mbox{for all created groups $j {\in\mathbb{N}}$},
\]
then either the cardinality of each (created) group tends to $\infty$,
a.s., or each of them will be eventually finite, a.s.
\end{longlist}

\end{theorem}
We show in Example~\ref{ex1} that for the collection of update
functions $f_{j}(n) = \e^{(j^{3} +n)}$, the cardinality of each group
remains finite, a.s. The third phase seems to arise only when for
fixed $n$, $j \to f_{j}(n)$ is an unbounded sequence.

The previous two theorems rely on a novel modification of a well-known
tool used in reinforced random walk processes, the Rubin construction,
which embeds $\operatorname{GAM}(\{f_{j}\},\{s_{n}\})$. We believe that such a
generalized Rubin construction as introduced in our paper could have
wider applicability to other preferential attachment models.

In the second part of the paper, we are going to estimate the
probability that a given group is the leading one. Our first result
concerns a \textit{reinforced urn model}. Consider an urn
with $k$ white balls and $ 1$ red ball and with reinforcement scheme
$f$. Then if we pick a ball at random, it is white with probability
$ f(k)/(f(k) + f(1))$, and red with probability $ f(1)/(f(k) + f(1))$.
Suppose that by the time of the $n$th extraction
we picked $j$ white balls and $n-j$ red ones. The probability to pick a
white ball becomes $ f(k+j)/(f(k+j) + f(n+1-j))$. We call the urn with
these initial conditions and dynamics a reinforced urn model with
parameters $k$ and $f$ [abbreviated $\operatorname{RUM}(k, f)$]. Denote by $\mathbb
{P}^{\ssup k}$ the probability measure referring to $\operatorname{RUM}(k, f)$. We
have the following estimate.
%
\begin{theorem}\label{estur} Fix any $k \ge1$ and consider a $\operatorname{RUM}(k,
f)$ with $\sum_{j=1}^{\infty} 1/\break f(j)<\infty$. We have
%
\begin{equation}
\label{bmsb} \qquad\mathbb{P}^{\ssup k} (\mbox{only a finite number of white balls
are picked} ) \le \frac12 \prod_{\ell=1}^{k-1}
\frac{f(\ell)F_{k}}{1
+ f(\ell) F_{k}},
\end{equation}
where $F_{k} \dff\sum_{j = k}^{\infty} 1/f(j)$.
\end{theorem}
The above theorem sheds deeper insight on the evolution of $\operatorname{RUM}(k,
f)$ and on Theorem~\ref{princres}(i): it shows that the leading side
in the beginning has a great probability to stay the dominant side. As
an example of the power of our bound, take $f(j)=j^2$. In this case, a
simple computation gives that
\begin{eqnarray*}
&&\mathbb{P}^{\ssup k} (\mbox{only a finite number of white balls are picked}
)\\
&&\qquad\le\frac{1}{2}\exp \biggl(-(k-1)+\frac{\pi}{2\sqrt{k+2}} \biggr).
\end{eqnarray*}
Hence for large initial weights $k$ the white has an overwhelming
chance to be the one with cardinality tending to infinity. The estimate
in \eqref{bmsb} improves Theorem~3.6.2 in~\cite{tzhou09}. Theorem \ref
{estur} should be also compared with Theorem~3 in~\cite{oliv}, which is
proved under the technical assumptions on the update function $f$
stated in Section 4 of that paper. Note also that the bound above is an
improvement on the upper bound which could be obtained in (\ref{bmsb})
by means of a similar reasoning to the one in Propositions~2.1 and~3.1 from~\cite{climic09}. The cause for this is that the
lower/upper bounds in~\cite{climic09} are rough for large initial
weights. This is one main reason why the methods there only work for
finite graphs and not also for infinite graphs. Our proof is based on
an embedding of $\operatorname{RUM}(k, f)$ into Brownian motion, and gives robust
estimates for all initial weights.

Next we turn again to $\operatorname{GAM}(f,\{s_{n}\})$. Suppose that $\sum_{j=1}^{\infty}1/f(j)<\infty$. Theorem~\ref{princres}(i) guarantees
the existence of a unique group whose cardinality goes to infinity. We
call this the leading group. Denote by $\mathit{Lead}$ the label of the
leading group. In other words, $\mathit{Lead} = j$ if and only if the
leading group is the $j$th one.
Our goal is to test if a given group, which has a certain advantage on
the others, is the leader. We start by giving an upper bound for the
tail of \textit{Lead}.

We give the following construction of $\operatorname{GAM}(f,p)$.
Suppose we have two sequences of random variables, $b_{n} $ and $
t(n)$, satisfying the following. The variables $b_{n}$ are i.i.d.
Bernoulli with mean $p$, while the variables $t(n)$ are described
recursively. We define $A^{*}_{1}(1)=1$, and $A^{*}_{i}(1)=0$ for all
$i \ge2$. Moreover, set $L^{*}_{1} =1$.
Denote by $\Fcal_{n}$ the $\sigma$-algebra generated by $ \{
(b_{i}, t(i) )$, with $i \le n \}$. Suppose we defined
$A^{*}_{i}(n)$, which is $\Fcal_{n-1}$-measurable. The random variable
$t(n)$ can be chosen to have the following distribution:
\[
\mathbb{P}\bigl(t(n) = k | \Fcal_{n-1}\bigr)= \frac{f_{k}(A^{*}_{k}(n))}{\sum_{s=1}^{L^{*}_{n}} f_{s}(A^{*}_{s}(n))}.
\]
Moreover, we can choose $t(n)$ to be independent of $b_{i}$ with $i\ge n+1$.
Denote by $L^{*}_{n} \df\max\{ j \ge1 \dvtx A^{*}_{j}(n) \ge1\}$.
We define
\[
A^{*}_{j}(n+1) \df \cases{ 0, & \quad $\forall j >
L^{*}_{n}+1,$\vspace*{2pt}
\cr
1, & \quad $\mbox{if
$b_{n} =1$},$\vspace*{2pt}
\cr
A^{*}_{j}(n) +
\1_{\{t(n) =j\}}, & \quad $\mbox{if $b_{n}=0$.}$ }
\]
Finally, let $L^{*}_{n+1} \df\max\{ j \ge1 \dvtx A^{*}_{j}(n+1) \ge
1\}$.
We have that $\{A^{*}_{i}(n),  i,n \in\mathbb{N}\}$ is distributed
like the process $\{A_{i}(n),  i,n \in\mathbb{N}\}$ described in
Section~\ref{MAIN}.
At time $n$, $b_{n}$ will determine if the new vertex is disconnected,
and $t(n)$ will determine to which of the existing vertices the new
arrived will adhere \textit{if} it is not disconnected. Notice that $t(n)$
is defined also in the case that $b_{n}=1$, that is, in the case that
the new vertex is disconnected.
We denote by $\xi_{1} = 0$ and $\xi_{i} \df\inf\{ n > \xi_{i-1}
\dvtx b_{n} =1\}$. In words, $\xi_{i} $ is the time when the $i$th
group is formed.
We say that the $i$th group is generated by the $u$th group if $t(\xi_{i})=u$; that is,
if we flipped the value of $b_{\xi_{i}}$ into~$0$,
then the new arrival would have joined the group $u$. In this case we
say that $u$ is the parent of $i$. Notice that there exists exactly one
parent for each integer different from one. We build a random tree
$\Gcal$, whose root is one, joining each integer to its parent. We say
that a vertex is at level $n$ if its distance from the root is $n$.
Denote by $g_{n}$ the vertices at level $n$. Let $G_{n} = \bigcup_{j \ge
n} g_{j}$. We have:
%
\begin{theorem}\label{testo1} Suppose that the assumptions of
Theorem~\ref{princres}\textup{(i)} hold. Then
\[
\mathbb{P}(\mathrm{Lead} \in G_{n}) \le\inf_{r,M \ge1}
\bigl[m^{n} \e^{-c_{n}(r,M) n} + r^{-n}+ C_{1}\exp
\{- M C_{2} \} \bigr],
\]
where the sequence $ c_{n}(r,M) \ti$ as $n \ti$, for fixed value of $r
\ge1$ and $M \ge1$, and $m, C_{1}, C_{2}>0$. The quantities $C_{1}$,
$C_{2}$ and $m$ are computable. The functions $c_{n}(r,M)$ are
computable for fixed values of $r$ and $M$.
\end{theorem}
The following result is a direct consequence of Theorems~\ref{testo1}
and~\ref{estur}.
%
\begin{cor}\label{rf} Suppose that the assumption of Theorem  \ref
{princres}\textup{(i)} holds. Then
\begin{eqnarray*}
\mathbb{P}(\mathrm{Lead} = 1)&\ge& 1 - \Biggl(\sum_{k=1}^{\infty}
\frac12 \prod_{\ell=1}^{k-1} \frac{f(\ell)F_{k}}{1 + f(\ell) F_{k}}
\Biggr)\\
 &&{}- \inf_{r,M \ge1} \bigl[m^{2} \e^{-2 c_{2}(r,M) } +
r^{-2} -C_{1}\exp \{- M C_{2} \} \bigr],
\end{eqnarray*}
where the quantities $c_{n}(r,M)$, $m$ and $C_{1}$ and $C_{2}$ are the
same as Theorem~\ref{testo1}.
\end{cor}

The rest of the paper is structured as follows: in Section \ref
{proof21} we introduce our generalized Rubin construction and give the
proof of Theorem~\ref{princres}(i). In Section~\ref{proof22} we give
the proof of Theorem~\ref{princres}(ii). In Section~\ref{proof23} we
prove our main result, Theorem~\ref{princres1}, and present an example
where the third phase occurs. In Section~\ref{prooft1} we introduce our
Brownian motion embedding and provide the proof of Theorem~\ref{estur}.
In Section~\ref{lead} we give the proofs of Theorem~\ref{testo1} and of
Corollary~\ref{rf}. Finally, in the \hyperref[app]{Appendix} we give a
brief introduction to the Rubin construction, as introduced in~\cite{BD1990}.

\section{\texorpdfstring{Proof of Theorem \protect\ref{princres}(i)}{Proof of Theorem 1.1(i)}}
\label{proof21}
We introduce a modified version of the Rubin construction which fits
our model. For a detailed explanation of the original Rubin
construction, see, for example,~\cite{climic09} and~\cite{BD1990}.

Fix a parameter $p<1$. We first focus on the case $s_{n}\equiv p<1$,
that is, $\operatorname{GAM}(f, p)$, and then we {extend} to the more general case
$s_{n} \le p$ using a coupling.
For any set $A \subset\R^{+}$, let
\[
A[n] = \inf\bigl\{ x \dvtx\# \bigl(A \cap[0,x]\bigr) \ge n+1\bigr\},
\]
where the infimum of an empty set is $\infty$. In words, $ A[n]$ is the
$n+1$th element of~$A$, ordered from the smallest to the largest. For
example, if $A = \{ 2, 8, 6, 9\}$, then $A[0]=2$ and $ A[1]= 6$, $A[5]
= \infty$. Notice also that for the example $A = \{ 1/j \dvtx j \ge1\}
$, is not possible to identify the $n+1$th element. In fact, in this
case, we have that $A[n] = 0$ for all $n \ge0$.

Notice that $A[n]$ is always a nondecreasing sequence, hence $\lim_{n
\ti} A[n]$ exists, possibly infinite.
For each $i \in\mathbb N, $ let $\{W_{n}^{\ssup i}, n\ge1\}$ be a
sequence of independent exponential$(1)$ random variables, with $ n \in
\mathbb{N}$. Moreover let $\{R^{\ssup i}_{n}, n \ge1\}$ be a sequence
of i.i.d. Bernoulli such that $ \mathbb{P}(R^{\ssup i}_{n} =1) =p$. We
are going to use these sequences to generate a $\operatorname{GAM}(f,p)$. The
Bernoullis will be used to create new groups, while the exponentials
play a central role in the allocation of new individuals into existing
groups. We are assuming that all the variables involved are independent
of each other. Set $N_{m}(1) = 1$, for all $m\ge1$. Then, for $n \ge
2$, let
\begin{eqnarray*}
N_{m}(n) &\df& 1+\#\bigl\{j\dvtx j \le n-1 \mbox{ such that }
R^{\ssup
m}_{j} =0\bigr\},
\\
\Xi_{1} &\df& \{0\}\cup \Biggl\{\sum_{i=1}^{n}
\frac{W_{i}^{\ssup
1}}{f (N_{1}(i) )} \dvtx n \ge1 \Biggr\} \subset\R^{+}.
\end{eqnarray*}
In words, for each $m \ge1$, the processes $\mathbf{N}_{m} \df \{
N_{m}(n), n \ge1\}$ are independent processes {with the property that
$N_{m}(n)-1$ are distributed like binomial with parameters $n-1$ and
$1-p$}, while
$\Xi_{1}$ is a random subset of $\R^{+}$ composed by $0$ and all the
partial sums of the sequence $ \{W_{i}^{\ssup1}/f (N_{1}(i)
),$ with $i \ge1 \}$. To each element $\Xi_{1}$ we associate a
corresponding Bernoulli as follows. Let $g_{1} \dvtx\Xi_{1} \to\{
0,1\}$ be a random function defined by $g_{1}(\Xi_{1}[n]) \df R^{\ssup1}_{n}$.
The elements in $\Xi_{1}$ {with} corresponding Bernoulli equal to one,
are used to generate new groups for $\operatorname{GAM}(f,p)$. The other ones will
\textit{potentially} belong to the first group and will be labeled one. We
will clarify the last sentence at the end of the construction.
Define
\[
\widetilde{\Xi}_{1} \df\{0\} \cup \Biggl\{\sum
_{i=1}^{n}\frac
{W_{i}^{\ssup1}}{f (N_{1}(i) )} \dvtx n \ge1 \mbox{ and }
R^{\ssup1}_{n} = 0 \Biggr\},
\]
that is, $\widetilde{\Xi}_{1}$ is composed of \{0\} and all the
points in $\Xi_{1} \setminus\{0\}$ {with} Bernoulli {equal to} $0$.
These are the points which do not generate other groups.\vadjust{\goodbreak}
We label the points in $\widetilde{\Xi}_{1}$ with $1$.
Set $\tau_{1} = 0 $ and define
\[
\tau_{2} \df\inf\bigl\{ n \ge1 \dvtx R^{\ssup1}_{n}
=1 \bigr\}.
\]
The random variable $\tau_{2}$ is the time when the second group is formed.
Given $\tau_{2}$, let
%
\[
\Xi_{2} \df\Xi_{1} \cup \Biggl\{ \Xi_{1}[
\tau_{2}] + \sum_{i=1}^{n}
\frac
{W_{i}^{\ssup2}}{f(N_{2}(i))}\dvtx n \ge0 \Biggr\}
\]
and
\[
\widetilde{\Xi}_{2} \df \Biggl\{\Xi_{1}[
\tau_{2}] + \sum_{i=1}^{n}
\frac
{W_{i}^{\ssup2}}{f (N_{2}(i) )} \dvtx\mbox{ either } n =0 \mbox { or both } n \ge1 \mbox{ and }
R^{\ssup2}_{n} = 0 \Biggr\}.
\]
We label the elements of $ \widetilde{\Xi}_{2} $ using 2.
Define the function $g_{2} \dvtx\Xi_{2} \to\{0,1\}$ as follows. If
$\Xi_{2}[n] = \Xi_{1}[j]$, for some $j \in\mathbb{N}$, then $g_{2}(\Xi_{2}[n]) = R^{\ssup1}_{j}$. The latter is well defined because all the
elements of $\Xi_{1}$ are a.s. distinct. If $\Xi_{2}[n] = (\Xi_{2}\setminus\Xi_{1})[j]$ for some $j \in\mathbb{N}$, then $g_{2}(\Xi_{2}[n]) = R^{\ssup2}_{j}$.
Notice that $\widetilde{\Xi}_{1}$ and $\widetilde{\Xi}_{2}$ are
disjoint, and their union is a proper subset of $\Xi_{2}$.
Let us describe in words the variables defined so far. The
reinforcement plays no role up to time $\tau_{2}$. The latter random
variable is geometrically distributed with mean $1/p$. At time $\tau_{2}$, the first group has cardinality $\tau_{2}$, because we count
also the point $0$, and a second group is formed. The random point $\Xi_{2}[\tau_{2}]$ is labeled 2, in fact it belongs to $\widetilde{\Xi
}_{2}$, and it is the smallest point belonging to this random set. The
next point on the line, that is, $\Xi_{2}[\tau_{2}+1]$ can have label
$1$, $2$ or no label at this stage. If the latter happens, we label
this point with~3. If it belongs to $\widetilde{\Xi}_{1}$,
respectively,
$\widetilde{\Xi}_{2}$, its label will be 1, respectively, 2. Notice that
by the definition of these sets, if $\Xi_{2}[\tau_{2}+1] \in\widetilde
{\Xi}_{1} \cup\widetilde{\Xi}_{2}$ then $g_{2} (\Xi_{2}[\tau_{2}+1] )$ must be equal to zero. On the other hand, in the case
that $g_{2} (\Xi_{2}[\tau_{2}+1] )=1$ then a new group is formed,
which is labeled $3$. The probability that this happens is $p$.
Next we want to compute the probability that $\Xi_{2}[\tau_{2}+1]$ has
label 1. We have the following equality:
%
\begin{eqnarray}
\label{eve}\qquad &&\bigl\{\Xi_{2}[\tau_{2}+1] \in\widetilde{
\Xi}_{1}\bigr\}
\nonumber
\\[-8pt]
\\[-8pt]
\nonumber
&&\qquad = \bigl\{\Xi_{1}[\tau_{2}+1]-
\Xi_{1}[\tau_{2}] < \widetilde{\Xi}_{2}[1] -
\Xi_{1}[\tau_{2}] \bigr\} \cap \bigl\{ g_{2}
\bigl(\Xi_{1}[\tau_{2}+1] \bigr) =0 \bigr\}.
\end{eqnarray}
Note that $\Xi_{2}[\tau_{2}] = \Xi_{1}[\tau_{2}]$.
Given $\tau_{2}$, the two events appearing on the right-hand side of
\eqref{eve} are independent, because the first one depends on the
exponentials while the second is determined by the Bernoullis. The
probability of the second event, conditionally on $\tau_{2}$, is $1-p$.
If the random variable $\Xi_{2}[\tau_{2}+1] $ was labeled~1, then it
would belong to $\widetilde{\Xi}_{1}$ and would be equal to
\[
\widetilde{\Xi}_{1}[\tau_{2}+1] = \sum
_{i=1}^{\tau_{2} +1} W_{i}^{\ssup
1}/f
\bigl(N_{1}(i) \bigr) = \Xi_{1}[\tau_{2}] +
W_{\tau_{2}+1}^{\ssup
1}/f \bigl(N_{1}(\tau_{2}+1)
\bigr).
\]
If $\Xi_{2}[\tau_{2}+1] $ was labeled 2, then it would belong to
$\widetilde{\Xi}_{2}$ and would be equal to $\widetilde{\Xi}_{2}[1] =
\Xi_{2}[\tau_{2}] + W_{1}^{\ssup2}/f (1 )$. Hence
\begin{eqnarray*} \Xi_{2}[\tau_{2}+1]&=&
\Xi_{2}[\tau_{2}] + \min \biggl(\frac{W_{\tau
_{2}+1}^{\ssup1}}{f (N_{1}(\tau_{2}+1) )},
\frac{W_{1}^{\ssup
2}}{f (1 )} \biggr)
\\
&=&\Xi_{2}[\tau_{2}] + \min \biggl(\frac{W_{\tau_{2}+1}^{\ssup1}}{f
(\tau_{2} )},
\frac{W_{1}^{\ssup2}}{f (1 )} \biggr),
\end{eqnarray*}
where we used $N_{1}(\tau_{2}+1) = \tau_{2}$. {This last} equality
comes from the fact that among~$R^{\ssup1}_{i}$, with $i \le\tau_{2}$, the only Bernoulli taking value one is $R^{\ssup1}_{\tau_{2}}$.
As $N_{1}(\tau_{2}+1)$ equals one plus the number of zeroes among the
first $\tau_{2}$ Bernoulli, it is equal to $\tau_{2}$. The first event
{on the right-hand side of} \eqref{eve} can be rewritten as
%
\begin{equation}
\label{eve1} \biggl\{\frac{W_{\tau_{2}+1}^{\ssup1}}{f (\tau_{2} )} < \frac
{W_{1}^{\ssup2}}{f (1 )} \biggr\}.
\end{equation}
Given $\tau_{2}$, the random variable $W_{\tau_{2}}^{\ssup1}/f (\tau_{2} ) $ is exponentially distributed with mean $1/f(\tau_{2})$. By
a simple integration, we can argue that the probability that, among two
independent exponentials, a given one is the smallest is equal to its
parameter divided by the sum of the parameters. Hence the probability
of the event in~\eqref{eve1}, conditionally on $\tau_{2}$, is $ f(\tau_{2})/ (f(\tau_{2}) + f(1)  )$. The probability of the event
described in~\eqref{eve}, {conditionally on $\tau_{2}$,} is
\[
(1-p) \frac{f(\tau_{2})}{f(\tau_{2}) + f(1) }.
\]
%
We infer that the {conditional} probability that $\Xi_{2}[\tau_{2}+1]
$ is labeled 2 is $ (1-p)f(1)/(f(\tau_{2}) + f(1))$. This is consistent
with what happens in $\operatorname{GAM}(f, p)$.

Suppose we defined $ (\tau_{2}, \Xi_{1}, \widetilde{\Xi}_{1}, g_{1},
\ldots,\tau_{m-1}, \Xi_{m-1}, \widetilde{\Xi}_{m-1}, g_{m-1} )$.
Define
\[
\tau_{m} \df\inf\bigl\{ n > \tau_{m-1} \dvtx
g_{m-1} \bigl(\Xi_{m-1}[n] \bigr) = 1\bigr\},
\]
that is, the time when the $m$th group is formed.
Given $\tau_{m}$ let
%
\begin{equation}
\label{xim} \Xi_{m} \df \Xi_{m-1} \cup \Biggl\{
\Xi_{m-1}[\tau_{m}] + \sum_{i=1}^{n}
\frac{W_{i}^{\ssup{m-1}}}{f(N_{m}(i))}\dvtx n\ge1 \Biggr\}
\end{equation}
and
%
\begin{eqnarray}
\label{ximtil}&& \widetilde{\Xi}_{m} \df \Biggl\{\Xi_{m-1}[
\tau_{m}] + \sum_{i=1}^{n}
\frac{W_{i}^{\ssup{m-1}}}{f(N_{m}(i))}\dvtx\mbox { either } n=0
\nonumber
\\[-8pt]
\\[-8pt]
\nonumber
&&\hspace*{72pt}\qquad{} \mbox{or both } n\ge1 \mbox{ and }
R^{\ssup m}_{n} = 0 \Biggr\}.
\end{eqnarray}
The elements of $ \widetilde{\Xi}_{m}$ are labeled $m$.
Moreover let $g_{m}$ be defined as follows. If there exists $j$ such
that $\Xi_{m}[n] =\Xi_{m-1}[j]$, then $g_{m} (\Xi_{m}[n] ) =
g_{m-1} (\Xi_{m-1}[j] )$. If $\Xi_{m}[n] = (\Xi_{m}\setminus\Xi_{m-1})[j] $ for some $j$, then set $ g_{m} (\Xi_{m}[n] ) =
R^{\ssup m}_{j}$.

Denote by $\Xi\df\bigcup_{s=1}^{\infty} \Xi_{s}$. Each point $x \in
\Xi$ belongs, a.s., to exactly one $\widetilde{\Xi}_{s}$ for some $s
\ge1$, that is, $ \Xi\df\bigcup_{s=1}^{\infty} \widetilde{\Xi}_{s}$.
In our construction, we label the point $x$ with~$s$ if and only if $x
\in\widetilde{\Xi}_{s}$. Define the random function $g \dvtx\Xi\to
\{0,1\}$ as follows. If $ \Xi[n] = \widetilde{\Xi}_{j}[s] $ for some
(a.s. unique) pair $(j,s) \in\mathbb{N}^{2}$, then $ g (\Xi[n] )
= R^{\ssup j}_{s}.$
Notice that $\Xi$ can be used to generate a generalized attachment
model, as follows. Denote by
\[
\widetilde{A}_{i}(n) = \bigl\{ j \dvtx j \le n, \Xi[j] \mbox{ has
label $i$}\bigr\}.
\]
{Then} $\{\widetilde{A}_{i}(n), \mbox{ with } i,n \ge1\}$ is
distributed like the process $\{A_{i}(n), \mbox{ with }\break i,n \ge1\}$
introduced in Section~\ref{MAIN}. To see this,
suppose that in the set $\{ \Xi[i], \mbox{ with } i \le n\}$ there are
exactly $\ell_{i}$ points labeled $i$, with $\sum_{i=1}^{m} \ell_{i} =n
$ for some $ m \in\mathbb{N} $ satisfying also $\ell_{i} \ge1$ for
all $i \in\{1,\ldots, m\}$. Given this, the probability that $\Xi
[n+1]$ is labeled $m+1$, that is, the probability that $g(\Xi[n+1])$
equals one,
is exactly $p$. Given that $\Xi[n+1]$ is not labeled $m+1$, then the
probability that it is labeled $j$, with $ j \le m$, is exactly
%
\begin{equation}
\label{permo} \frac{f(\ell_{j})}{\sum_{i=1}^{m} f(\ell_{i})},
\end{equation}
where we used the memoryless property of the exponential random
variables. In fact, using this property, given that $\Xi[n+1]$ is not
labeled $m+1$, the random variable $\Xi[n+1]- \Xi[n]$ is distributed
like the minimum of $m$ exponentials with parameters $f(\ell_{s})$, for
$1 \le s \le m$. The probability that the $j$th exponential is the
minimum is given exactly by \eqref{permo} through a simple integration.
Summarizing, given that in the set $\{ \Xi[i],$ with $ i \le n\}$ there
are exactly $\ell_{i}$ points labeled $i$, with $\sum_{i=1}^{m} \ell_{i} =n$ and $\sum_{i=1}^{m-1} \ell_{i} < n$ for some $m \in\mathbb
{N}$, the probability that $\Xi[n+1]$ is labeled $j$, with $ j \le m$, is
\[
(1-p) \frac{f(\ell_{j})}{\sum_{i=1}^{m} f(\ell_{i})}.
\]

Define
%
\begin{equation}
\label{xstar} x^{*}_{m} \df\Xi_{m-1}[
\tau_{m}] + \sum_{i=1}^{\infty}
\frac
{W_{i}^{\ssup{m}}}{f(N_{m}(i))},
\end{equation}
and for any integer $j \ge1$, set
%
\begin{equation}
\label{star} \Xi^{*}_{j} \df \Biggl\{ \Xi[
\tau_{j}] + \sum_{s=1}^{n}
W^{\ssup
j}_{s}/f \bigl(N_{j}(s) \bigr) \dvtx n \ge0
\Biggr\}.
\end{equation}
In the next result, we prove that $ x^{*}_{m} $ is a.s. finite, for any
$m\ge1$. This, together with~\eqref{xim} and \eqref{ximtil},\vadjust{\goodbreak} implies
that $ x^{*}_{m}$ is an accumulation point for $\Xi_{m} $ and
$\widetilde{\Xi}_{m}$.
We say that a vertex $u$ is generated by $j$ if $\Xi[\tau_{u}] \in\Xi^{*}_{j}$. Notice that each vertex (different from 1) is generated by
exactly one other vertex.
Our proof of Lemma~\ref{finit} relies on the construction of a random
tree $\mathcal{T}$, built by connecting each vertex to its parent.
Notice that this random tree shares the same distribution with $\Gcal$,
introduced before Theorem~\ref{testo1}. Suppose that $\tau_{u} = t$. If
we switched $g( \Xi[t])$ from~1 to 0, we would have that $\Xi[t]$ would
have been a point of $\widetilde{\Xi}_{j}$, and hence it would have had
label $j$. Fix $j,n \in\mathbb{N}$. Notice that even if the Bernoulli
associated to the point $\Xi^{*}_{j}[n]$ equals~1, this point might not
be able to generate a child in $\mathcal{T}$ using the exponentials and
Bernoulli that have been defined so far. This is the case if $\#
(\Xi\cap[0, \Xi^{*}_{j}[n]] ) = \infty$, {when} infinitely many
vertices {have already been} generated by the time we reach $\Xi^{*}_{j}[n]$ and all the $(W^{\ssup i}_{n}, R^{\ssup i}_{n})$ have
already been used. This is going to be an important point in the proof
of Lemma~\ref{genac}.

\begin{Lemma}\label{finit}
The random variables $x^{*}_{m}$, with $m \ge1$, are almost surely finite.
\end{Lemma}
\begin{pf}Fix $m \ge1$. Set $Z_{m}(0)=0 $, and let $Z_{m}(i) = \inf\{
n \dvtx N_{m}(n)=i\}$. Then $Z_{m}(i) - Z_{m}(i-1)$, with $i \ge1$
are $\operatorname{geometric}(1-p)$ and are independent of the $W_{i}^{\ssup m}$, with
$i \ge1$. If $Z_{m}(i) \le k < Z_{m}(i+1)$, then $f (N_{m}(k) )
= f(i)$. Hence
\begin{eqnarray*} x^{*}_{m} &=& \Xi[
\tau_{m}] + \sum_{i=1}^{\infty} \sum
_{j=
Z_{m}(i)}^{Z_{m}(i+1) -1}\frac{W_{j}^{\ssup{m}}}{f(N_{m}(j))}
\\
&=& \Xi[\tau_{m}] + \sum_{i=1}^{\infty}
\frac1{f(i)} \sum_{j=
Z_{m}(i)}^{Z_{m}(i+1) -1}W_{j}^{\ssup{m}}.
\end{eqnarray*}
As the series in the latter expression is composed by nonnegative
random variables, it is a.s. finite if its mean is finite. Its mean is exactly
%
\begin{equation}
\label{finf} \frac1{1-p} \sum_{i=1}^{\infty}
\frac1{f(i)}<\infty.
\end{equation}
To see this, notice that $Z_{m}(i),  i \ge1$, is independent of
$W^{\ssup m}_{j}$, $ j \ge1$, which implies
\[
\mathbb{E}\Biggl[ \sum_{j= Z_{m}(i)}^{Z_{m}(i+1) -1}W_{j}^{\ssup{m}}
\Biggr] = \frac1{1-p}.
\]
Moreover, we have that $\Xi[\tau_{m}]$ is stochastically smaller than
%
\begin{equation}
\label{sumf} \frac{1}{\min_{i} f(i)} \sum_{s=1}^{\tau_{m}}W^{\ssup s}_{1}.
\end{equation}
This is because $\Xi[n] - \Xi[n-1]$ is stochastically smaller than an
exponential random variable whose mean is smaller than $1/(\min_{i}
f(i))$. Moreover, the random variable $\tau_{m}$ is negative binomial
with parameters $m$ and $p$. This can be checked by induction; in fact,
$\tau_{1}$ is geometrically distributed with mean $1/p$. Suppose {this}
is true for $\tau_{m-1}$. Then we have to wait for {an}
independent $\operatorname{geometric}(p)$ to create the next group.
Combining this fact with \eqref{sumf} we have that $\Xi[\tau_{m}]<\infty
$ a.s.
This, together with \eqref{finf} implies the lemma.
\end{pf}

In the next result we establish the link between the behavior of the
generalized attachment model and the quantity $ \inf_{i} x^{*}_{i}$.

\begin{Lemma}\label{genac} The infimum $ \inf_{i} x^{*}_{i}$, is a.s.
attained, that is, it is actually a minimum. The minimizer is a.s.
unique. Moreover
%
\begin{equation}
\label{info} \lim_{n \ti} \Xi[n] = \inf_{i}
x^{*}_{i} \qquad\mbox{a.s.}
\end{equation}
\end{Lemma}
\begin{pf}
We select a random subtree of $\mathcal{T}$, denoted by $\mathcal
{T}_{1}$, as follows.
The root of this tree is 1 (i.e., {it} is identified with the first
group). Given that a vertex $j$ belongs to $\mathcal{T}_{1}$, its
offspring will be those vertices $u$ such that
%
\begin{equation}
\label{require} \Xi[\tau_{u}] \in\Xi^{*}_{j}\quad
\mbox{and}\quad {x^{*}_{u}} < x^{*}_{j}.
\end{equation}
Recall that $\Xi[\tau_{u}] = \Xi_{u-1}[\tau_{u}]= \Xi_{u}^{*}[0].$
We are going to prove the following statement:
%
\begin{equation}\label{fqaf}
\begin{tabular}{p{280pt}@{}}
For any fixed $M$, only finitely many of the vertices $u$ of $
\mathcal{T}_{1}$ satisfy $\widetilde{\Xi}_{u}[1] < M$.
\end{tabular}
\end{equation}
Before we prove \eqref{fqaf} we argue that this statement would imply
the lemma. {We need only consider the vertices of $\mathcal{T}_{1}$. In
fact, if $j$ is not a vertex of $\mathcal{T}_{1}$, then there exists a
vertex $u$ such that $x^{*}_{j} > x^{*}_{u}$. Hence $x^{*}_{j} \neq\inf_{i} x^{*}_{i}$.}

If \eqref{fqaf} holds, then for any $M$ there are only finitely many
vertices $u$ in $\mathcal{T}_{1}$ such that $x^{*}_{u}<M$. Hence, as
each $x^{*}_{u}$ is a.s. finite, we have that
$\inf_{i} x^{*}_{i}$ is actually a minimum. Next we prove that the
minimizer is a.s. unique. To prove this last statement, we prove that
for each $i>j$, we have that $x^{*}_{i}$ and $x^{*}_{j}$ are a.s.
different. To see this, notice that $x^{*}_{i} - \Xi[\tau_{i}]$ only
depends on $\{W^{\ssup i}_{n}, R^{\ssup i}_{n}, $ with $ n\ge1\}$.
Hence $ x^{*}_{i} - \Xi[\tau_{i}]$ is independent of $ x^{*}_{j} - \Xi
[\tau_{i}]$ which is determined by a disjoint collection of
exponentials and Bernoullis. The probability that $ x^{*}_{i} - \Xi[\tau_{i}]$ and $ x^{*}_{j} - \Xi[\tau_{i}]$ are equal is 0, as they are
continuous independent random variables. This is exactly the
probability that $x^{*}_{i} = x^{*}_{j}$. As the set of $x^{*}_{i},  i
\ge1$, is countable,
$x^{*}_{i}$ are all, a.s., distinct.

Next we show that \eqref{fqaf} implies \eqref{info}. As already
mentioned, the sequence $\Xi[n]$ is a.s. nondecreasing, that is, $\Xi
[n+1] \ge\Xi[n]$, a.s. Hence\break $\lim_{n \ti} \Xi[n]$ a.s. exists.\vadjust{\goodbreak} Notice
that for each $i$, $x^{*}_{i}$ is the limit of an increasing sequence
taking values in $\Xi$.
To see this, notice that
\[
\Xi[\tau_{i}] + \sum_{j=1}^{n}
\frac{W_{j}^{\ssup{i}}}{f(N_{i}(j))} < x^{*}_{i}\qquad \forall n \ge1,
\]
by the definition of $ x^{*}_{i}$. Hence infinitely many points labeled
$i$ are smaller than $ x^{*}_{i}$, yielding $\# (\Xi\cap[0,
x^{*}_{i}] ) = \infty$.
This implies that $\lim_{n \ti} \Xi[n] \le x_{i}^{*}$ for each $i\ge
1$, that is, $\lim_{n \ti} \Xi[n] \le\inf_{i} x_{i}^{*}$.
Now we turn to the proof of the other inequality which implies~\eqref
{info}. Fix $\eps>0$. It is sufficient to prove that~\eqref{fqaf} implies
%
\begin{equation}
\label{acm1} \#\Bigl({(u,j) \dvtx\widetilde{\Xi}_{u}[j]} \le
\inf_{i} x^{*}_{i} - \eps \Bigr) < \infty.
\end{equation}
In fact, if \eqref{acm1} holds, only finitely many $u$ satisfy $\#
(\widetilde{\Xi}_{u} \cap [0, \inf_{i} x^{*}_{i} - \eps ] )$
{is infinite}. Denote the set of labels of these groups by $B$. For
each element $u$ of the finite set $B$, there are only finitely many
points of $\Xi^{*}_{u}$ which are smaller than $\inf_{i} x^{*}_{i} -
\eps$, for otherwise we would have $x^{*}_{u} \le \inf_{i} x^{*}_{i} -
\eps$ which would yield a contradiction. Hence, for each element $u$
of $B$, the set $\Xi^{*}_{u}\cap [\inf_{i} x^{*}_{i} - \eps ] $
is finite. For each $j \notin B$ we have that there exists a $u \in B$
such that $\widetilde{\Xi}_{j}[0] \in\Xi^{*}_{u}$. This implies that
\[
\Xi\cap \Bigl[0, \inf_{i} x^{*}_{i} - \eps
\Bigr] = \bigcup_{u \in B} \Xi^{*}_{u}
\cap\Bigl[\inf_{i} x^{*}_{i} - \eps \Bigr],
\]
and the latter is a finite set.
Notice that $ \inf_{i} x^{*}_{i} - \eps$ could take a negative value,
but this is not a problem for our reasoning, as {then} the set
appearing in \eqref{acm1} would be empty, and there would be nothing to prove.
Next, we prove that \eqref{fqaf} implies \eqref{acm1}. Fix a vertex $u$
of $\mathcal{T}_{1}$. Denote by $u(0)=1, u(1), u(2), \ldots, u(n)$ the
ancestors of $u$ in $\mathcal{T}_{1}$, that is, the vertices lying on
the unique self-avoiding path connecting $u$ to the root $1$. Notice
that we do not consider $u$ ancestor of itself. If $u$ satisfies
$\widetilde{\Xi}_{u}[1] \le \inf_{i} x^{*}_{i} - \eps$, then $u$
belongs to $\mathcal{T}_{1}$.
In fact, $\widetilde{\Xi}_{u}[1] \ge\Xi^{*}_{u(i)}[1]$ while $ \inf_{i} x^{*}_{i} - \eps< x^{*}_{u(i-1)}$, for all $i \le n+1$. Hence, $
\Xi^{*}_{u(i)}[1] \le x^{*}_{u(i-1)},$ where we set $u(n+1) = u$. As we
are assuming that~\eqref{fqaf} holds, the random tree $\mathcal{T}_{1}$
has only finitely many vertices $j $ satisfying $\widetilde{\Xi
}_{j}[1]>M$. As $\inf_{i} x^{*}_{i} <\infty$, a.s., we have that~\eqref
{acm1} holds.

Next, we are going to prove \eqref{fqaf}.
For any vertex $j$ in $\mathcal{T}_{1}$, denote by $\sigma_{j}$ the
number of its offspring. Notice that the $\sigma_{j}$ are neither
independent nor identically distributed, and
$\mathcal{T}_{1}$ is not a Galton--Watson tree. To see this, fix $j,n
\ge1$. If there is an infinite number of elements of $\Xi$ to the left
of $\Xi^{*}_{j}[n]$, that is,
%
\begin{equation}
\label{inft1} \# \bigl(\Xi\cap \bigl[0, \Xi^{*}_{j}[n]
\bigr] \bigr) = \infty,
\end{equation}
then already infinitely many groups have been created. Hence $\Xi^{*}_{u}[0] < \Xi^{*}_{j}[n]$ for all $u \in\mathbb{N}$. This implies
that $\Xi^{*}_{j}[n]$ cannot generate any new group\vadjust{\goodbreak} in $\mathcal
{T}_{1}$ using the exponentials and Bernoullis defined so far, because
they have already been used. To overcome this problem, we create a new
tree, larger than $\mathcal{T}_{1}$, by introducing new random
variables which allow also the observations $\Xi^{*}_{j}[n]$ satisfying
\eqref{inft1} to create a new group. To this end, we should attach to
each of these observations a new sequence of
independent exponentials and independent Bernoullis. For example, if
$\Xi^{*}_{j}[n]$ satisfies \eqref{inft1}, and the associated Bernoulli
equals one, {a} new group, that we label $\nu$, is created (notice
that we cannot use any of the integers as a label, because they are
already all taken). In this case, we set $\widetilde{\Xi}_{\nu}[0] = \Xi^{*}_{j}[n]$.
We denote the associated sequence of i.i.d. exponentials
with mean 1 by $W^{\ssup\nu}_{n}$, and let $R^{\ssup\nu}_{j}$ be the
Bernoulli associated to group $\nu$. {We} define $\Xi^{*}_{\nu}$ and
$\widetilde{\Xi}_{\nu}$ using these random variables, as we did in
\eqref{star} and \eqref{ximtil}. If the group $\nu$ satisfies the
second requirement in \eqref{require}, that is,
${x^{*}_{\nu} - \Xi^{*}_{\nu}[0]}< x^{*}_{j} - \Xi^{*}_{\nu}[0]$, then
$\nu$ belongs to the new tree $\mathcal{T}_{2}$ that we are going to
define. But then we {would} have to allow {that} $\nu$ {is} able to
generate groups as well, in the same fashion. This approach would
require that we introduce new sequences of exponentials and Bernoullis,
and the notation would be quite awkward. Hence we prefer a different
approach. Before we proceed in a formal description of $\mathcal
{T}_{2}$ , notice that for this tree the number of offspring per vertex
are independent and identically distributed. {In fact, } $ {x^{*}_{\nu}
- \Xi^{*}_{\nu}[0]}$ is independent of $x^{*}_{j} - \Xi^{*}_{\nu}[0]$ {
as determined by disjoint sets of exponentials and Bernoulli}.
Moreover, analyzing the event $\{{x^{*}_{\nu} - \Xi^{*}_{\nu}[0]} <
x^{*}_{j} - \Xi^{*}_{\nu}[0]\}$, one can easily argue that it does not
depend on the exponentials and Bernoullis attached to vertices
different from $u$ and $j$. Summarizing, the number of offspring of $j$
in this new tree {depends} only:
\begin{itemize}
\item on the exponentials attached to $j$, with the exception of
$W^{\ssup j}_{1}$;
\item on the Bernoullis attached to $j$;
\item and on $W^{\ssup\nu}_{1}$, if $ \nu= \Xi^{*}_{j}[n]$ for some
$n$ and $R^{\ssup j}_{n} =1$.
\end{itemize}
This implies that the number of offspring per vertex are i.i.d.

Now we are ready to give a formal construction of $\mathcal{T}_{2}$.
Suppose that to each $x \in\Xi^{*}_{1}$ we associate an extra {sequence of exponential random variables $\Theta^{\ssup x}_{i}$, with
parameter 1, and an independent copy, say $N_{x}(i)$, of $N_{1}(i)$,
with $i \ge1$ } Let
\[
\eta_{1} \df\#\Biggl\{x \in\Xi^{*}_{1} \dvtx{
\sum_{i=1}^{\infty}\Theta^{\ssup x}_{i}/f
\bigl(N_{x}(i)\bigr)} < x^{*}_{1} - x \mbox{ and
} { g_1}(x)=1\Biggr\}.
\]
The previous random variable counts also the $n$ satysfying $\# (\Xi
\cap[0, \Xi^{*}_{1}[n]] ) = \infty$, hence $\eta_{1}$ is
stochastically larger than $\sigma_{i}$ for any $i$. Then the
Galton--Watson tree $\mathcal{T}_{2}$ whose offspring distribution is
the same as the one of $\eta_{1}$ is stochastically larger than
$\mathcal{T}_{1}$. We assume that $\mathcal{T}_{2}$ is built on the
same probability space of $\mathcal{T}_{1}$. In other words, we can
assume, and we will, that $\mathcal{T}_{1}$ is a subtree of $\mathcal{T}_{2}$.
Next we prove that the
average number of offspring is bounded by a finite constant $m$.
Define
%
\begin{equation}
\label{mb1} \Omega_{k,j} \df \Biggl\{{\sum
_{i=1}^{\infty} W^{\ssup
k}_{i}/f
\bigl(N_{k}(i)\bigr)}< \sum_{s=k+1}^{\infty}W^{\ssup j}_{s}/f
\bigl(N_{j}(s) \bigr) \Biggr\}.
\end{equation}
Notice that we should have used different $\operatorname{exponentials}(1$) instead of
$W^{\ssup k}_{1}$, but the two share the same distribution and are
independent of the right-hand side, and this notation is easier to
handle. Of course, we are allowed to do that because we are interested
only in estimating the probability of this event.

We have that
%
\begin{equation}
\label{om1} \mathbb{E}[\eta_{j}] \le\mathbb{E}\Biggl[\sum
_{k=1}^{\infty} \1_{\Omega
_{k,j}}\Biggr] = \sum
_{k=1}^{\infty} \mathbb{P}(\Omega_{k,j}).
\end{equation}
In order to prove \eqref{om1}, notice that {on} the left-hand side we
count the number of elements in $\Xi^{*}_{u}$ {with} Bernoulli equal
to 1, and {which} satisfy an extra condition. The right-hand side
counts only those vertices which satisfy the extra condition.
Hence we only need to prove that $\mathbb{P}(\Omega_{k,j})$ is
summable. Notice that $\mathbb{P}(\Omega_{k,j})$ is independent of $j$.

Recall that
\[
\theta^2_{k} =1\Big/ \Biggl(\sum_{s=k+1}^{\infty}
2/(1-p)^2f^2(s) \Biggr).
\]
Denote by $Y_{k} = \sum_{s=k+1}^{\infty}W^{\ssup j}_{s}/f
(N_{j}(s) )$, and {$Z = \sum_{i=1}^{\infty} W^{\ssup
k}_{i}/f(N_{k}(i))$}. We have
%
\begin{eqnarray}
\label{oho} %
 E\bigl[\e^{\theta_{k} Y_{k}}\bigr] &=& \prod
_{s=k+1}^{\infty} \sum_{j=1}^{\infty
}
\biggl(\frac{f(s)}{f(s) - \theta_{k}} \biggr)^{j} p^{j-1}(1-p)
\nonumber\\
&=& \prod_{s=k+1}^{\infty} \biggl(
\frac{f(s)}{f(s) - \theta_{k}} \biggr) (1-p) \frac{1}{1-  (pf(s)/(f(s) - \theta_{k}) )}
\\
&=& \prod_{s=k+1}^{\infty}\frac{ f(s)(1-p)}{f(s)(1-p) -
\theta_{k}}\nonumber
\end{eqnarray}
and
\begin{eqnarray*}
\mathbb{E}\bigl[\e^{- \theta_{k} Z}\bigr] &= &\prod_{s=1}^{\infty}
\frac
{f(s)(1-p)}{f(s)(1-p) + \theta_{k}}\\
&=& \prod_{s=1}^k
\frac
{f(s)(1-p)}{f(s)(1-p) + \theta_{k}} \prod_{s \ge k+1} \frac
{f(s)(1-p)}{f(s)(1-p) + \theta_{k}}.
\end{eqnarray*}
Hence,
%
\begin{eqnarray}\label{we1} %
&&\mathbb{E}\bigl[\e^{- \theta_{k} Z}\bigr] E
\bigl[\e^{\theta_{k} Y_{k}}\bigr]\nonumber\\
&&\qquad\le \prod_{s=1}^k
\frac{f(s)(1-p)}{f(s)(1-p) + \theta_{k}} \prod_{s \ge k+1} \frac{f^2(s)(1-p)^2}{f^2(s)(1-p)^2 - \theta^2_{k}}
\nonumber\\
&&\qquad= \prod_{s=1}^k \frac{1}{1 +  (\theta_{k}/(f(s)(1-p)) )} \prod
_{s=k+1}^{\infty}\frac{1}{1 -  (\theta^2_{k}/(f^2(s)(1-p)^2) )}
\\
&&\qquad\le\prod_{s=1}^k \frac{1}{1 +  (\theta_{k}/(f(s)(1-p)) )} \exp
\Biggl\{ C \theta_{k} \sum_{s=k+1}^{\infty} 1/f^2(s)\Biggr\}\nonumber
\\
&&\qquad\le {\mathrm{const}} \prod_{s=1}^k
\frac{1}{1 +  (\theta_{k}/(f(s)(1-p))
)}, \nonumber
\end{eqnarray}
where {for the first inequality} we used that {$\theta_{k}/((1-p)^2f^2(s)) \le1/2$ for $s \ge k$ and our choice of $\theta_{k}$, and
the inequality $1-x \ge\e^{-Cx}$, for $x \in(0,1/2)$ for a
proper choice of $ C$. Using the assumptions in Theorem \ref
{princres}(i) and (\ref{we1}), we have that}
%
\begin{equation}
\label{mo} \mathbb{E}[\eta_{j}] \le \mathbb{E}\Biggl[\sum
_{k=1}^{\infty} \1_{\Omega
_{k,j}}\Biggr] \le{\sum
_{k=1}^{\infty} \mathbb{E}\bigl[\e^{- \theta_{k} Z}\bigr] E
\bigl[\e^{\theta_{k} Y_{k}}\bigr]} \df m<\infty,
\end{equation}
where, for the finitess of $m$ we used the second assumption in
Theorem~\ref{princres}(i), and the fact that $\Gamma$ is a finite set.

For each vertex $u $ in $ \mathcal{T}_{2}$, recall that we denote by
$\Xi^{*}_{u}[0]$ the time when this vertex was generated and by $\Xi^{*}_{u}[n] = \Xi^{*}_{u}[0]+ \sum_{j=1}^{n} W^{\ssup
u}_{j}/f(N_{u}(j))$. This is consistent with our definition given in
\eqref{star}, but now it is defined for indices which are not
necessarily integers.
Next we prove that each vertex $u$ at level $n+1$ has a probability to
satisfy $\widetilde{\Xi}_{u}[1]< M$ which decreases faster than $\e^{-cn}$ for any $c>0$.
For any vertex $u \in\mathcal{T}_{2}$ we denote by $|u|$ its distance
from the root of the tree. Recall that the set of vertices at distance
$k$ from the root is called level $k$. Fix a large parameter~$M$. A
vertex $u$ of $ \mathcal{T}_{2}$ is \textit{good} if the element which
generates $u$ is smaller than~$M$.
A~path is a (possibly finite) sequence of vertices $u(i), i \ge1$,
such that $u(i+1)$ is generated by $u(i)$. We say that a path connects
vertex $a$ to level $n+1$ if the first element of the path is $a$ and
the last {lies} at level $n+1$. We build the following random path
$\mathbf{u}$. We start from $1=u(0)$ and if this vertex has at least
one offspring in $\mathcal{T}_{2}$, we choose one at random assigning
the same probability to each {offspring}. We denote its label as
$u(1)$. If $u(1)$ has at least one offspring, we choose one of them at
random and denote its label by $u(2)$. We follow this procedure until
we either reach level $n+1$ or find a vertex with no offspring. The
event $\{$the path~$\mathbf{u}$ connects $1$ to a vertex at level $n+1\}$
equals the\vadjust{\goodbreak} event that each of the $u(i)$ has at least one offspring. Hence
\[
\mbox{$ \{$the path $\mathbf{u}$ connects $1$ to a vertex at level $n+1 \}$}=
\bigcap_{i=0}^{n} \{\eta_{u(i)}
\ge1 \}.
\]
%
Notice that each event $  \{\eta_{u(i)} \ge1 \}$ is independent
of $\Xi^{*}_{u(i-1)}[1]$ and is independent of each $W^{\ssup\ell
}_{k}$ with $\ell< u(i-1)$, and $ k \ge1$. Moreover the events $
\{\eta_{u(i)} \ge1 \}$ are independent.
Define
\[
\Psi(n,k) \df \Biggl\{ \sum_{i=1}^{n}
\1_{\{ \eta_{u(i)} \le k\}} \ge0.5 n \Biggr\}.
\]
Fix $k \ge1$. We have
%
\begin{eqnarray}
\label{nz1} %
&&\mathbb{P}\bigl(u(n+1) \mbox{is good }
| \mbox{the path $\mathbf{u}$ connects $1$ to a vertex at level $n+1$}\bigr)
\nonumber\hspace*{-35pt}\\
&&\qquad\le \mathbb{P} \Biggl(\sum_{i=1}^{n}
\Xi^{*}_{u(i)}[1] - \Xi^{*}_{u(i)}[0] \le
M {|} \mbox{the path $\mathbf{u}$ connects $1$}\nonumber\hspace*{-35pt}\\
&&\hspace*{135pt}\qquad\mbox{to a vertex at level $n+1$}
\Biggr)
\nonumber\hspace*{-35pt}\\
&&\qquad= \mathbb{P} \Biggl(\sum_{i=1}^{n}
\Xi^{*}_{u(i)}[1] - \Xi^{*}_{u(i)}[0] \le
M \Big| \bigcap_{i=1}^{n} \{
\eta_{u(i)} \ge1 \} \Biggr)\hspace*{-35pt}
\\
&&\qquad\le\mathbb{P} \Biggl(\sum_{i=1}^{n} \bigl(
\Xi^{*}_{u(i)}[1] - \Xi^{*}_{u(i)}[0]
\bigr)\1_{\{\eta_{u(i)} \le k \}} \le M \Big|\bigcap_{i=1}^{n}
\{\eta_{u(i)} \ge1 \} \Biggr)
\nonumber\hspace*{-35pt}\\
&&\qquad\le\mathbb{P} \Biggl(\sum_{i=1}^{n} \bigl(
\Xi^{*}_{u(i)}[1] - \Xi^{*}_{u(i)}[0]
\bigr)\1_{\{\eta_{u(i)} \le k \}} \le M \Big| \bigcap_{i=1}^{n}
\{\eta_{u(i)} \ge1 \} \cap\Psi(n,k) \Biggr)
\nonumber\hspace*{-35pt}\\
&&\qquad\quad{} + \mathbb{P} \Biggl(\Psi^{c}(n,k) \Big| \bigcap
_{i=1}^{n} \{\eta_{u(i)} \ge1 \} \Biggr).\nonumber\hspace*{-35pt}
\end{eqnarray}
In the last step we used that for any triplet of events $A$, $B$, $C$
we have
\[
\mathbb{P}(A | B) \le\mathbb{P}(A | B\cap C) + \mathbb{P}\bigl(C^{c}
| B\bigr).
\]
Next, we bound the last probability in \eqref{nz1},
%
\begin{eqnarray}
\label{nz0} %
&&\mathbb{P} \Biggl(\Gamma^{c}(n,k)
\Big| \bigcap_{i=1}^{n} \{\eta_{u(i)}
\ge1 \} \Biggr)
\nonumber\\
&&\qquad= \mathbb{P} \Biggl(\sum_{i=1}^{n} \1\{
\eta_{u(i)} \le k \} \le0.5 n \Big| \bigcap_{i=1}^{n}
\{ \eta_{u(i)} \ge1\} \Biggr)
\\
&&\qquad= \mathbb{P} \Biggl(\sum_{i=1}^{n} \1\{
\eta_{u(i)} > k \} \ge0.5 n \Big| \bigcap_{i=1}^{n}
\{ \eta_{u(i)} \ge1\} \Biggr).\nonumber
\end{eqnarray}
Let $\xi_{u(i)}$, $i \le n$, be i.i.d. random variables taking values
in $\mathbb{N}$, with distribution
\[
\mathbb{P}(\xi_{u(i)} \ge k) = \mathbb{P}(\eta_{u(i)} \ge k |
\eta_{u(i)} \ge1)= \frac{\mathbb{P}(\eta_{u(i)} \ge k)}{\mathbb{P}(\eta_{u(i)} \ge1)} \df q_{k} \qquad\mbox{for } k
\ge1.
\]
The sequence $q_{k}$ is independent of $u(i)$ because the random
variables $\eta_{i}$ are i.i.d.
Moreover, as the $\eta_{u(i)}$ are independent, $\sum_{i=1}^{n} \1\{
\eta_{u(i)} > k \}$ is, conditionally on $\bigcap_{i=1}^{n} \{ \eta_{u(i)}
\ge1\}$, binomially distributed with mean $ n q_{k}$.
If $X$ is a binomial with parameters $(n,q)$, then
%
\begin{equation}\qquad
\mathbb{P}(X \ge0.5 n) \le\exp \biggl\{ - \biggl( \frac{1}{2q} \ln\biggl(
\frac
{1}{2q}\biggr) + \frac1{2(1-q)} \ln\frac1{2(1-q)} \biggr) n \biggr\},
\end{equation}
by a simple exponential bound; see, e.g.,~\cite{dembo88} pages 27 and 35.

Fix $ r >1$. We can choose $K^{*}_r$ such that $q_{K^{*}_r} < 1/2$ and
\begin{eqnarray*}
&&\mathbb{P} \Biggl(\sum_{i=1}^{n} \1\bigl\{
\xi_{u(i)} > K^{*}_r\bigr\} \ge0.5 n \Biggr)\\
&&\qquad \le
\exp \biggl\{ - \biggl( \frac{1}{2q_{K^{*}_r}} \ln\biggl(\frac
{1}{2q_{K^{*}_r}}\biggr) +
\frac1{2(1-q_{K^{*}_r})} \ln\frac1{2(1-q_{K^{*}_r})} \biggr) n \biggr
\} \\
&&\qquad\le \frac1{(r m)^{n}},
\end{eqnarray*}
where $m$ has been defined in \eqref{mo}.
We can choose such $ K^{*}_{r} $ because\break $\lim_{k \ti}(1/(2q_{k})) \ln
(1/(2q_{k})) = \infty$. Notice that for any $k \ge K^{*}_r$, we have $
q_{k} \le q_{K^{*}_r}< 1/2$. Moreover if $k\ge K^{*}_r$, then
\begin{eqnarray*}
&&\exp \biggl\{ - \biggl( \frac{1}{2q_{k}} \ln\biggl(\frac{1}{2q_{k}}\biggr) +
\frac1{2(1-q_{k})} \ln\frac1{2(1-q_{k})} \biggr) n \biggr\}\\
&&\qquad
\le \frac1{(r m)^{n}}.
\end{eqnarray*}
This fact is due to the monotonicity of $q_{k}$ and the convexity of
the function $2x \ln(2x) + 2(1-x) \ln2(1-x)$, for $ x \in(0,1)$, and
the fact that this function attains its minimum at $1/2$.
Next, let {$ (e_{i} )$} be a sequence of i.i.d. random
variables with distribution
\begin{eqnarray*}
&&\mathbb{P}(e_{i} \le x ) = \mathbb{P} \Biggl(W^{\ssup2}_{1}/f(1)
\le x {|} W^{\ssup2}_{1}/f(1)\\
&&\hspace*{47pt}\qquad \le\sum
_{t=K^{*}_r+1}^{\infty} W^{\ssup1}_{t}/f
\bigl(N_{1}(t)\bigr) - \sum_{j =2}^{\infty}
W^{\ssup2}_{j} f\bigl(N_{2}(j)\bigr) \Biggr).
\end{eqnarray*}
In words, $e_{i}$ is distributed like an exponential with mean $
1/f(1)$ conditioned to be smaller than {an independent
quantity}.
We claim that the first probability in the last equation of \eqref{nz1}
is smaller or equal to
%
\begin{equation}
\label{nz2} \mathbb{P} \Biggl(\sum_{j=1}^{\lfloor 0.5 n \rfloor}
e_{i}\le M \Biggr).
\end{equation}
To see this, notice that by a simple exchangeability argument we have that
\begin{eqnarray*}
&& \mathbb{P} \Biggl(\sum_{i=1}^{n}
\bigl(\Xi^{*}_{u(i)}[1] - \Xi^{*}_{u(i)}[0]
\bigr)\1_{\{\eta_{u(i)} \le K^{*}_r \}} \le M {\Big|} \bigcap_{i=1}^{n}
\{\eta_{u(i)} \ge1\} \cap\Psi\bigl(n,K^{*}_r
\bigr) \Biggr)
\\
&&\qquad= \mathbb{P} \Biggl(\sum_{i=1}^{n} \bigl(
\Xi^{*}_{u(i)}[1] - \Xi^{*}_{u(i)}[0]
\bigr)\1_{\{\eta_{u(i)} \le K^{*}_r \}} \le M {\Big|} \bigcap_{i=1}^{ \lfloor 0.5 n \rfloor }
\bigl\{1 \le\eta_{u(i)} \le K^{*}_r\bigr\}
\Biggr)
\\
&&\qquad\le\mathbb{P} \Biggl(\sum_{i=1}^{\lfloor 0.5n \rfloor }
\bigl(\Xi^{*}_{u(i)}[1] - \Xi^{*}_{u(i)}[0]
\bigr)\1_{\{\eta_{u(i)} \le K^{*}_r \}} \le M {\Big|} \bigcap_{i=1}^{\lfloor 0.5 n \rfloor }
\bigl\{1 \le\eta_{u(i)} \le K^{*}_r\bigr\}
\Biggr).
\end{eqnarray*}
Again, notice that the events $ \{1 \le\eta_{u(i)} \le K^{*}_r\}$,
with $ i \le n+1$, are independent. Given $ \{1 \le\eta_{u(i)} \le
K^{*}_r\}$, the random variable $\Xi^{*}_{u(i)}[1] - \Xi^{*}_{u(i)}[0]
$ is stochastically larger than $e_{i}$, as $\sum_{s=k}^{\infty}
W^{\ssup1}_{s}/f(N_{1}(s))$ is a.s. decreasing in $k$. This proves the
relationship between \eqref{nz2} and the first probability in the last
equation of \eqref{nz1}.
Next a simple exponential bound, which uses the fact that $e_{i}$ are
independent, yields
\begin{eqnarray*} \mathbb{P} \Biggl(\sum_{j=1}^{\lfloor 0.5 n \rfloor }
e_{i}\le M \Biggr)&=& \mathbb {P} \Biggl(\frac1{0.5 n}\sum
_{j=1}^{\lfloor 0.5 n1 \rfloor } e_{i}\le \frac1{0.5 M} \Biggr)
\\
&=& \mathbb{P} \Biggl(\exp \Biggl\{ - \theta\frac1{0.5 n}\sum
_{j=1}^{\lfloor 0.5 n \rfloor} e_{i} \Biggr\}\ge \exp\biggl\{
- \theta\frac1{0.5 M} \biggr\} \Biggr)
\\
&\le&\exp\bigl\{- c_{n}(r,M) n\bigr\},
 \end{eqnarray*}
where $c_{n}(r,M) \ti$ as $n \ti$. For each $n$, $c_{n}(r,M)$ is the
Fenchel--Legendre transform (i.e., we minimize the exponent on $\theta
$) of $e_{i}$ in the point $\frac1{0.5 M} $.
Hence, the number of good vertices in $\mathcal{T}_{1}$ at level $n$ is
smaller or equal to
%
\begin{equation}
\label{cmlf} m^{n} \biggl( \exp \bigl\{-c_{n}(r,M) n \bigr\}
+ \frac{1}{(r m)^{n}} \biggr) .
\end{equation}
Hence only finitely many vertices in $\mathcal{T}_{2}$ are good. This
implies that only finitely many vertices in $\mathcal{T}_{1}$ are good,
and this, in turn, implies \eqref{fqaf}.
\end{pf}

\begin{pf*}{Proof of Theorem~\ref{princres}(\textup{i})} First suppose that
$s_{n} \equiv p<1$. The minimum of $\inf_{i} x^{*}_{i}$ is a.s. unique,
and we denote\vadjust{\goodbreak} it by $J^{*}$. By Lemma~\ref{genac} $\lim_{n \ti} \Xi[n]
= x^{*}_{J^{*}}$, hence the cardinality of group $J^{*}$ tends to
infinity, while the cardinality of {each of} the {other} groups {is} finite.

Now we reason for general $s_{n} \le p$, using a simple coupling. Let
$\{S_{i}, i\ge1\}$ be a sequence of independent Bernoullis with
$\mathbb{P}(S_{i} =1) = s_{i}/p = 1 - \mathbb{P}(S_{i} =0)$. We use
these random variables to relabel the points in $\Xi$ as follows.
If $S_{1} =0$, then we set $\Theta_{1} = \Xi_{2}[\tau_{2}] \cup\Xi
\setminus\widetilde{\Xi}_{2}$. If $S_{1} =1$, then $\Theta_{1} = \Xi$.
Define $\widetilde{\tau}_{3} \df\inf\{ n > \tau_{2} \dvtx g(\Theta_{1}[n])=1\}$. Suppose we have defined $\Theta_{m-1}$ and $\widetilde
{\tau}_{i}$, for $i \le m$. On the event $\{\sum_{i=1}^{m-1}S_{i}=k\},$
if $S_{m} =0$, respectively, $S_{m} =1$, set $\Theta_{m} = \Xi_{k+1}[\widetilde{\tau}_{k+1}] \cup\Theta_{m-1} \setminus\widetilde
{\Xi}_{k+1}$, respectively, $\Theta_{m} = \Theta_{m-1}$. We set
$\widetilde{\tau}_{m+1} \df \inf\{ n > \widetilde{\tau}_{m} \dvtx
g(\Theta_{m}[n])=1\}$. Let $\Theta= \bigcap_{i} \Theta_{i}$. The
process $\Theta[n]$, with $n \ge1$ is a $\operatorname{GAM}(f, \{s_{n}\}$). Let
$\kappa(n) = \sum_{j=1}^{n} S_{j}$. Denote by $h(i) = \inf\{ n \dvtx
\kappa(n) =i\}$. This implies that
the $i$th group in $\Theta$ is the $h(i)$th group in $\Xi$. Let $U_{i}
\df\widetilde{\Xi}_{h(i)}$, and
%
\begin{equation}
\label{usp2} u^{*}_{i} \df x^{*}_{h(i)}.
\end{equation}
This implies that $\inf_{j} \{ u^{*}_{j} \dvtx j \ge1\}$ is actually
a minimum and has a unique minimizer. Following the same reasoning
given in the previous paragraph, we conclude that the only group whose
cardinality grows to infinity is $K^{*}_r$.
\end{pf*}

\begin{pf*}{Proof of Corollary~\ref{findegree}}
We first assume that $s_{n} \equiv p$. For any $i$, denote by $E(i)$
the set of groups which are generated by $i$. In virtue of \eqref{mo},
we have that
%
\begin{equation}
\label{ny1} \qquad V(u)\df\bigl\{\widetilde{\Xi}_{i}[1] <
x^{*}_{u} \mbox{ for only finitely many $i \in E(u)$}\bigr
\} \qquad\mbox{holds a.s.}
\end{equation}
Notice that for $u$ which is not a vertex of $\mathcal{T}_{1}$ we have
that $\widetilde{\Xi}_{u}[1] > \inf_{i}x^{*}_{i}= \lim_{n \ti} \Xi[n]$.
Hence, we {do not have} to consider such $u$.
Recall the definition of $G_{N}$ given before Theorem~\ref{testo1}. As
for each {$N$}, there are only finitely many good vertices in $\mathcal
{T}_{1}$, and {we get}
$\lim_{N \ti}\mathbb{P}(\mathrm{Lead} \in G_{N}) =0$. Combining the
latter limit with \eqref{ny1} we have that
\begin{eqnarray*}
&&\mathbb{P}\Bigl(\lim_{n \ti}
A_{u}(n)>1 \mbox{ for only finitely many $u$}\Bigr)
\\
&&\qquad= \mathbb{P}\Bigl(\widetilde{\Xi}_{u}[1] < \inf_{i}
x^{*}_{i} \mbox{ for only finitely many $u$}\Bigr)
\\
&&\qquad= \lim_{N \ti} \mathbb{P}\Bigl(\Bigl\{\widetilde{\Xi}_{u}[1]
< \inf_{i} x^{*}_{i} \mbox{ for only finitely
many $u$}\Bigr\}\cap\{\mathrm{Lead} \notin G_{N}\}\Bigr)
\\
&&\qquad\ge\lim_{N \ti} \mathbb{P} \biggl(\bigcap_{u \in\mathcal{T}_{1} \dvtx
u \notin G_{N}}
V(u) \cap\{\mathrm{Lead} \notin G_{N}\} \biggr)
\\
&&\qquad= \lim_{N \ti} \mathbb{P}(\mathrm{Lead} \notin G_{N}) =1.
\end{eqnarray*}
For the general case $s_{n} \le p$, apply the same coupling we used at
the end of the previous proof.
\end{pf*}

\section{\texorpdfstring{Proof of Theorem \protect\ref{princres}(ii)}{Proof of Theorem 1.1(ii)}}\label{proof22}
We first deal {with} the case $s_{n} = p$. Repeat the construction
given in the proof of Theorem~\ref{princres}(i), under the hypothesis
of Theorem~\ref{princres}(ii). Recall the definition of $\Xi_{u}^{*}$, $\widetilde{\Xi}_{u}$ and $x^{*}_{u}$. Recall also the
definition of $\mathcal{T}$. The random variables $x^{*}_{u}$ , for $u
\ge1$, are a.s. infinite, because the infinite sum of independent
exponentials is finite if and only if its mean is finite.
We prove {next} that for any fixed $M>0$,
%
\begin{equation}
\label{liM} \liminf_{u \ti} \widetilde{\Xi}_{u}[1] > M\qquad
\mbox{a.s.}
\end{equation}
Fix a vertex $u_{n}$ of $\mathcal{T}$, and denote by $u_{i}$, with $i
\le n-1$ its ancestors: that is, $\Xi_{u_{j}}[\tau_{u_{j}}] \in\Xi^{*}_{u_{j-1}}$, for all $j \le n$. Then $\Xi_{u_{n}}[1]$ is
stochastically larger than a sum of $n-1$ i.i.d. exponentials with
parameter $f(1)$. Hence $\lim_{n \ti} \Xi_{u_{n}}[1] = \infty$, a.s.
Now notice that $\Xi^{*}_{s}\cap[0,M]$ is a.s. finite for each $s \ge
1$. Hence, as $u$ grows to infinity, the number of its ancestors grows
to infinity,
proving \eqref{liM}.
Since it is easy to adapt the above reasoning to the case $s_{n} \le
p$, we will leave this task to the reader.

\section{\texorpdfstring{Proof of Theorem \protect\ref{princres1}}{Proof of Theorem 1.5}}
\label{proof23}
\mbox{}
\begin{pf*}{Proof of Theorem~\ref{princres1}(\textup{i})} We first analyze the
case $s_{n} \equiv p$. We build a similar construction {as the one}
given in the proof of Theorem~\ref{princres}, {the only difference
being that we } place $f_{i}$ instead of $f$. We leave to the reader to
check that this construction embeds our $\operatorname{GAM}(\{f_{j}\}, p)$. In this setting,
\[
x^{*}_{i} \df\Xi_{i}[\tau_{i}] +
\sum_{j=1}^{\infty}\frac{W^{\ssup
i}_{j}}{f_{i} (N_{i}(j) )}.
\]
Notice that $ x^{*}_{i}$ is a.s. finite if and only if $ \sum_{s=1}^{\infty} 1/f_{i}(s)$ is finite. Hence, we do not exclude that $
x^{*}_{i} = \infty$, a.s., but we know that
%
\begin{equation}
\label{acpoin} \mbox{there exists at least one $j$ for which
$x^{*}_{j}<\infty$, a.s.}
\end{equation}

Denote by $y^{*}$ the smallest accumulation point of $\Xi$. This
minimum accumulation point exists because the set of accumulation
points of $\Xi$ is closed, and the set $\Xi$ is a subset of $\R^{+}$.
Moreover $y^{*}$ is a.s. finite because of \eqref{acpoin}.
If $y^{*} < x^{*}_{i}$ for all $i \ge1$, then
%
\begin{equation}
\label{acp} \lim_{n \ti}\Xi[n] <x^{*}_{i}\qquad
\mbox{for all $i \ge1$.}
\end{equation}
We need to prove \eqref{acp} only for the case $x^{*}_{i} < \infty$,
because for the other cases the result is implied by the fact that
$y^{*}$ is an accumulation point which is a.s. finite and $\# \Xi\cap
[0, y^{*}+ \eps] = \infty$. Assume that $x^{*}_{i}$ is a.s. finite and
notice that for fixed $i$, as $y^{*} < x^{*}_{i}$, then $\delta_{i} \df
(x^{*}_{i} +y^{*})/2 < x^{*}_{i}$. As $y^{*}$ is an accumulation point
for $\Xi$, then $ \#\Xi\cap[0, \delta_{i}]$ is a.s. infinite. In words\vadjust{\goodbreak}
there are infinitely many points of $\Xi$ smaller than $\delta_{i}$. Hence
%
\begin{equation}
\label{acp1} \lim_{n \ti} \Xi[n] \le\delta_{i} <
x^{*}_{i}.
\end{equation}
The inequality in \eqref{acp1} holds for each $i$, yielding \eqref{acp}.
Moreover, \eqref{acp1} implies that each group will end up having
finite cardinality. This is because $\#\widetilde{\Xi}_{i} \cap[0,
\delta_{i}]$, as $\delta_{i}$ is strictly less than $x^{*}_{i}$ which
is the only accumulation point for $\widetilde{\Xi}_{i}$. The latter
statement is a direct consequence of the definitions of $x^{*}_{i}$
and~$\widetilde{\Xi}_{i}$.\looseness=-1

On the other hand, if $y^{*} = x^{*}_{i}$ for some $i$, then using
again that all the $x^{*}_{j}$ which are finite are also a.s. distinct,
we have that
$\lim_{n \ti}\Xi[n] = x^{*}_{i}.$ To prove the latter inequality,
suppose it is not true, that is, $\lim_{n \ti}\Xi[n] < x^{*}_{i}.$ Then
there would be an accumulation point smaller than $y^{*}$, which would
yield a contradiction.

Next we analyze the general case, that is, $s_{n} \le p$, for some
$p<1$ and all $n\ge1$.
The problem here is that the reinforcement function is group dependent.
In the special case $s_{n} \equiv p$ we had that the first point
labeled $i$ was $\Xi^{*}_{i}[0]$. We need to translate the points
labeled $i$ in the new construction for the general case. Denote by
$v(i)$ the time when the $i$th group is created and denote by $\Upsilon_{i-1}$ the union of the points labeled $j$, with $j \le i-1$. We have
that the first point labeled $i$ is exactly $\Upsilon_{i-1}[v(i)].$ Set
\begin{eqnarray}\label{usp3}
U^{*}_{i} &\df& \bigl\{ \Xi^{*}_{i}[n]
- \Xi^{*}_{i}[0]+ \Upsilon_{i-1}\bigl[v(i)\bigr]
\dvtx n \ge0 \bigr\},
\nonumber
\\[-8pt]
\\[-8pt]
\nonumber
 u^{*}_{i} &\df& x^{*}_{i}
- \Xi^{*}_{i}[0]+ \Upsilon_{i-1}\bigl[v(i)\bigr].
\end{eqnarray}
Hence $\Upsilon_{i} = \Upsilon_{i-1} \cup U^{*}_{i}.$ Moreover, let
$\Upsilon= \bigcup_{i=1}^{\infty}\Upsilon_{i}$. It is easy to check that
$ \Upsilon$ embeds $\operatorname{GAM}(\{f_{i}\}, \{s_{n}\})$.
We prove the theorem on the event $\{u^{*}_{j} $ is a.s. finite for at
least one created group $j\}$. Repeating the argument we gave for the
case $s_{n} \equiv p$, we see that either all the groups remains finite
or there exists exactly one dominating the others.
\end{pf*}

\begin{pf*}{Proof of Theorem~\ref{princres1}(\textup{ii})} First assume that
$s_{n} \equiv p$, for some $p<1$. Under the assumptions of this part of
the theorem, we have that each $x^{*}_{j} = \infty$, a.s. Hence $\inf_{j} x^{*}_{j} = \infty$. By our construction, either $\lim_{n \ti} \Xi
[n] = \infty$, in which case the cardinality of each group is a.s.
diverging to infinity, or $\lim_{n \ti} \Xi[n]= \Gamma< \infty$, a.s.,
in which case $\# \widetilde{\Xi}_{u} \cap[0, \Gamma]<\infty$, a.s. In
words, in the latter case, the cardinality of each group will
eventually remain finite, for otherwise $x^{*}_{j} \le\Gamma<\infty$
for some $j$, and this would give a contradiction.

For general $s_{n} \le p$, we have that $u^{*}_{i} = \infty$, where the
$u^{*}_{i}$ are the random variables defined in \eqref{usp3}. Reasoning
as in the previous paragraph we get the result for this more general case.
\end{pf*}

\subsection{An example when the third phase occurs}
\label{ex1}
Next we show an example where a third phase occurs, that is,
%
\begin{equation}
\label{thirdphase} \lim_{n \ti} A_{i}(n) < \infty \qquad\mbox{a.s.
for each $i \ge1$}.\vadjust{\goodbreak}
\end{equation}
In this example we pick $f_{j}(n) = \e^{(j^{3} +n)}$ and $s_{n} \equiv
p\in(0,1)$.
Notice that $\tau_{i+1} - \tau_{i}$, with $i \ge1$, is an i.i.d.
sequence of geometrically distributed random variables, with mean
$1/p$. Hence, by a standard exponential bound, we have
\[
\mathbb{P} \bigl(\tau_{n} > \bigl((1/p)- \eps \bigr) n \bigr) =
\mathbb{P} \Biggl(\sum_{i=1}^{n} (
\tau_{i} - \tau_{i-1}) > \bigl((1/p)- \eps \bigr) n \Biggr)
\le\e^{-C n}.
\]
This implies that
%
\begin{equation}
\label{sum1} \sum_{n=1}^{\infty}\mathbb{P}\bigl(
\tau_{n} > n^{2}\bigr) < \infty.
\end{equation}
Next, we use this fact to prove that
%
\begin{equation}
\label{intmd1} \mbox{for each $j\ge1$ there exists an $s >j$ such that
$x^{*}_{s} < x^{*}_{j}$. }
\end{equation}
The latter implies that $\inf_{j} x^{*}_{j}$ is not attained. As this
infimum is an accumulation point for $\Xi$, this would imply that the
smallest accumulation point of $\Xi$ is smaller than $x^{*}_{j}$, for
all $j \ge1$. Hence, \eqref{thirdphase} would hold.

Next we turn to the proof of \eqref{intmd1}. Fix $j \in\mathbb{N}$.
As $\Xi_{j} \subset\Xi$, we have that $\Xi[\tau_{u}] \le\Xi_{j}[\tau_{u}]$. Hence
%
\begin{eqnarray}
\label{seq} %
&&\mathbb{P} \bigl(x^{*}_{u}
> x^{*}_{j} \bigr)
\nonumber\hspace*{-35pt}\\
&&\quad=\mathbb{P} \Biggl(\Xi[\tau_{u}] + \sum
_{\ell=1}^{\infty} W^{\ssup
u}_{\ell}/f_{u}
\bigl(N_{u}(\ell)\bigr)> x^{*}_{j} \Biggr)
\nonumber\hspace*{-35pt}
\\[-8pt]
\\[-8pt]
\nonumber
&&\quad\le\mathbb{P} \Biggl(\Xi_{j}[\tau_{u}] + \sum
_{\ell=1}^{\infty} W^{\ssup
u}_{\ell}/f_{u}
\bigl(N_{u}(\ell)\bigr)> x^{*}_{j} \Biggr)\hspace*{-35pt}
\\
&&\quad\le\mathbb{P} \Biggl(\sum_{\ell=1}^{\infty}
W^{\ssup u}_{\ell
}/f_{u}\bigl(N_{u}(\ell)
\bigr) > \sum_{\ell=u^{2}}^{\infty} W^{\ssup j}_{\ell
}/f_{j}
\bigl(N_{j}(\ell)\bigr) {|} \tau_{u} < u^{2}
\Biggr)+ \mathbb{P} \bigl(\tau_{u} \ge u^{2} \bigr).\nonumber\hspace*{-35pt}
\end{eqnarray}
The last inequality in \eqref{seq} is justified as follows. For any
pair of events $A$ and $B$ we have that
\[
\mathbb{P}(A) = \mathbb{P}(A | B)\mathbb{P}(B) + \mathbb{P}\bigl(A\cap
B^{c}\bigr) \le \mathbb{P}(A | B) + \mathbb{P}\bigl(B^{c}
\bigr).
\]

Notice that $\{\tau_{u} < u^{2}\}$ is measurable with respect to the
$\sigma$-algebra
\[
\sigma \bigl\{R^{\ssup t}_{i} \dvtx t < u \mbox{ and } i <
u^{2} \bigr\}.
\]
In words, if we know the first $u^{2}-1$ observations of each
$\widetilde{\Xi}_{t}$, with $ t<u$, and the associated Bernoullis, we
know if the event
$\{\tau_{u} < u^{2}\}$ holds. Hence the latter event is independent of
the pair
\[
\Biggl(\sum_{\ell=1}^{\infty} W^{\ssup u}_{\ell}/f_{u}
\bigl(N_{u}(\ell)\bigr), \sum_{\ell=u^{2}}^{\infty}
W^{\ssup j}_{\ell}/f_{j}\bigl(N_{j}(\ell)
\bigr)\Biggr).
\]
Hence the last expression in \eqref{seq} equals
\[
\mathbb{P} \Biggl(\sum_{\ell=1}^{\infty}
W^{\ssup u}_{\ell
}/f_{u}\bigl(N_{u}(\ell)
\bigr) > \sum_{\ell=u^{2}}^{\infty} W^{\ssup j}_{\ell
}/f_{j}
\bigl(N_{j}(\ell)\bigr) \Biggr)+ \mathbb{P} \bigl(\tau_{u}
> u^{2} \bigr).
\]
The last expression is summable. To see this, in virtue of \eqref
{sum1}, we just need to prove that the first term is summable. Then our
argument follows from an application of the first Borel--Cantelli
lemma. In fact, the summability implies that $\{x^{*}_{u} < x^{*}_{j}\}
$ for infinitely many $u$. Set $\gamma_{u,j}= (1/j) \e^{-j^{3}
-u^{2}}$, and recall that $j$ is fixed.
For any pair of random variables $X$ and $Y$ and any constant $a$, we
have that
\begin{eqnarray*}
\mathbb{P}(X>Y) &=& \mathbb{P}(X>Y, X>a) +
\mathbb{P}(X>Y, X<a)\\
 &\le& \mathbb{P}(X>Y, X>a) + \mathbb{P}(Y<a)
\\
&\le&\mathbb{P}(X>a) + \mathbb{P}(Y<a).
\end{eqnarray*}
We apply this fact to obtain
\begin{eqnarray*} &&\mathbb{P} \Biggl(\sum
_{\ell=1}^{\infty} W^{\ssup u}_{\ell
}/f_{u}
\bigl(N_{u}(\ell)\bigr) > \sum_{\ell=u^{2}}^{\infty}
W^{\ssup j}_{\ell
}/f_{j}\bigl(N_{j}(\ell)
\bigr) \Biggr)
\\
&&\qquad \le\mathbb{P} \Biggl(\sum_{\ell=1}^{\infty}
W^{\ssup u}_{\ell
}/f_{u}\bigl(N_{u}(\ell)
\bigr) > \gamma_{u,j} \Biggr) + \mathbb{P} \Biggl(\sum
_{\ell
=u^{2}}^{\infty} W^{\ssup j}_{\ell}/f_{j}
\bigl(N_{j}(\ell)\bigr) < \gamma_{u,j} \Biggr).
 \end{eqnarray*}
Notice
\[
\sum_{n=1}^{\infty} 1/ f_{u}(n) =
\sum_{n=1}^{\infty} \e^{-u^{3} -n} =
\e^{-u^{3}} \sum_{n=1}^{\infty}
\e^{-n}= C_{1} \e^{- u^{3}},
\]
while, by a similar reasoning, $\sum_{n=u^{2}} 1/ f_{j}(n) \sim C_{2} \e^{- u^{2} - j^{3}}.$
Notice {that in virtue of Markov's inequality, we have }
\[
\mathbb{P} \Biggl(\sum_{\ell=1}^{\infty}
W^{\ssup u}_{\ell
}/f_{u}\bigl(N_{u}(\ell)
\bigr) > \gamma_{u,j} \Biggr) {\le C_{1} \e^{-u^{3}}/
\gamma_{u,j} = C_{1} j \exp\bigl\{-u^{3} +
j^{3}+u^{2}\bigr\}}
\]
and the right-hand side is summable in $u$ for fixed $j$.
In a similar way, using Chebyshev's inequality after applying the
function $\e^{\theta x}$ to both sides and choosing $\theta= (1-p)^{2}
\e^{u^{2}+j^{3}}$, we obtain
\[
\mathbb{P} \Biggl(\sum_{\ell=u^{2}}^{\infty}
W^{\ssup j}_{\ell
}/f_{j}\bigl(N_{j}(\ell)
\bigr) < \gamma_{u,j} \Biggr) \le\exp\bigl\{ - \bigl(\e^{ u^{2} +
j^{3}}
\bigr) ( C_{2} \e^{- u^{2} - j^{3}} - \gamma_{u,j}\bigr\}.
\]
The last expression is summable in $u$, because, for fixed $j$, $C_{2}
\e^{- u^{2} - j^{3}}$ is larger than $ \gamma_{u,j}$ for all
sufficiently large $u$.

\section{Brownian motion embedding}\label{prooft1}
Suppose that the positive function $f$ satisfies the condition $\sum_{j=1}^{\infty}1/f(j)<\infty$. Consider an urn
with $k$ white balls and $1$ red one. We pick a ball at random, and it
is white with probability
$ f(k)/(f(k) + f(1))$. Suppose that by the time of the $n$th extraction,
we picked $j$ white balls and $n-j$ red ones. The probability to pick a
white ball at the next stage becomes $ f(k+j)/(f(k+j) + f(n+1-j))$.
Let
\[
D \df\{\mbox{only a finite number of white balls are picked}\}.
\]
Denote by $\mathbb{P}^{\ssup k}$ the probability measures referring to
the urn with initial conditions and dynamics described above.

Let $F \df \sum_{j=1}^{\infty}1/f(j)$ and recall that $F_{k} \df\sum_{j=k}^{\infty}1/f(j)$ .
Let the process $\mathbf{B} := \{B_t,   t\ge
0\}$ be a standard Brownian motion, which starts from the point $F-
F_{k}=\sum_{i=1}^{k-1} 1/f(i)$. Denote by $\mathbb{Q}^{\ssup k}$ the
measure associated with this Brownian motion. We use this process to
generate the urn sequence described at the beginning of this section,
as follows. 
Set $m_{0} =0$ and let
%
\begin{equation}
\label{ennezero} m_{1} \df\inf \bigl\{ t \ge0\dvtx B_t-
{B_{0}} \mbox{ hits either } 1/f(k) \mbox{ or } - 1/f(1) \bigr\}.
\end{equation}

If $ B_{m_{1}} - B_{m_{0}} >0 $, then set $z_{1} =1$; otherwise set $z_{1}=0$.

Suppose we defined $m_{n}$ and $z_{1}, z_{2}, \ldots, z_{n}$. Set $\phi
(n) = \sum_{i=1}^{n}z_{i}$. On the event $\phi(n) = s$, we define
\[
m_{n+1} = \inf \biggl\{ t \ge m_{n}\dvtx B_t -
B_{m_{n}} \mbox{ hits either } \frac{1}{f(s+k)} \mbox{ or } -
\frac{1}{f(n-s+1) } \biggr\}.
\]
Set
\[
z_{n+1} \df \cases{
1, & \quad $\mbox{if $
B_{m_{n+1}} - B_{m_{n}}= 1/f(k+s)$},$
\vspace*{2pt}\cr
0, & \quad$ \mbox{if $ B_{m_{n+1}} - B_{m_{n}}= - 1/f(n-s +1) $}.$}
\]

By the ruin problem for Brownian motion, we have that
\begin{eqnarray*}
\mathbf{P} \bigl( z_{n+1}=0 | \phi(n) = s \bigr) &=& \frac{
1/f(s+k)}{ (1/f(s+k)) + (1/f(n-s+1))} \\
&=&
\frac{f(n-s+1)}{f(s+k) + f(n-s+1)},
\end{eqnarray*}
which is exactly the urn transition probability.\vadjust{\goodbreak}


In this way we embedded the urn into Brownian motion. In fact,
the process $\phi(n)$, with $ n \ge1$,
is distributed like the number of white balls withdrawn from the urn
associated to the reinforcement scheme described at the beginning of
this section.
Notice that
%
\begin{equation}
\label{punto} B_{m_n}= \sum_{j=1}^{k+ \phi(n)}
\bigl(1/f(j)\bigr) - \sum_{s=1}^{n - \phi(n)}
\bigl(1/f(s)\bigr)\qquad \mbox{with } n \ge0.
\end{equation}
Define
%
\begin{equation}
\label{enne} S \df\lim_{n \to\infty} m_n.
\end{equation}
This limit exists because the sequence of stopping times $\{m_{n}\}$ is
increasing. For this reason $S$ is itself a stopping time.
Define
\begin{eqnarray*}
D_{1} &\df& \Biggl\{ \exists n \ge1 \dvtx B_S = \sum
_{j=1}^{k+n} \bigl(1/f(j)\bigr) - \sum
_{j=1}^{\infty} \bigl(1/f(j)\bigr) \Biggr\} = \{
B_{S} <0 \} ,
\\
D_{2} &\df& \Biggl\{ \exists n\ge1 \dvtx B_S = \sum
_{j=1}^{\infty} \bigl(1/f(j)\bigr) - \sum
_{j=1}^{n} \bigl(1/f(j)\bigr) \Biggr\} = \{
B_{S} >0 \}.
\end{eqnarray*}
Moreover, in virtue of Theorem~\ref{hrubin} we have that exactly one
of the collection of events $ \{ \{z_{i}=0\}, i\ge1 \}$ and $
\{ \{z_{i}=1\}, i \ge1 \}$ holds finitely many times, a.s. This
implies that the event
$ D_{1} \cup D_{2} $
holds $\mathbb{Q}^{\ssup k}$-a.s. By our embedding, we have that
\[
\mathbb{Q}^{\ssup k} (D_{1}) = \mathbb{P}^{\ssup k} (D),
\]
where $D$ was defined at the beginning of this section.

\begin{pf*}{Proof of Theorem~\ref{estur}} In order to prove our result
we only need to prove the following:
\[
\mathbb{Q}^{\ssup k}(D_1) \le\frac12 \prod
_{s=1}^{k-1} \frac
{f(s)F_{k}}{1 + f(s) F_{k}}.
\]
Let
%
\begin{equation}
\label{dstu0} T \df\inf \biggl\{ n \ge1 \dvtx\phi(n) = \frac{n-k}2 \biggr
\}.
\end{equation}
This stopping time can be infinite with positive probability.
Notice that on $\{T<\infty\}$, by \eqref{punto}, we have that the urn
generated by the Brownian motion contains, at time $T$, an equal number
of white and red balls, and
$ B_{m_{T}} = 0.$ Viceversa, if we let
\[
H \df\inf\{ t\ge0 \dvtx B_{t} = 0\},
\]
then we have that
%
\begin{equation}
\label{tie1} \{H<S\} = \{T < \infty\}.
\end{equation}
To prove \eqref{tie1}, notice that for $k \in\mathbb{N}$, with $ k>0$,
the random sequence
\[
n \to \sum_{j=1}^{k+\phi(n)}\frac1{f(j)}- \sum
_{j=1}^{n - \phi
(n)}\frac1{f(j)}
\]
cannot switch sign without becoming $0$. So if $B_{m_{j}}>0$ and $
B_{m_{t}}<0$, for some $j < t$, then there exists an $s$, with $j < s <
t$, such that
$B_{m_{s}} = 0$. In this case, by time $s$ we have a tie. We use this
fact throughout the proof.

Recall that under $\mathbb{Q}^{\ssup k}$ the Brownian motion $\mathbf
B$ starts from $F- F_{k}$.
For $j \le k$, let
\[
H_{j} \df\inf\{ t \ge0 \dvtx B_{t} = F_{j+1}
-F_{k}\}.
\]
Notice that $F_{j+1} - F_{k}\ge0$ for $j \le k$. Moreover, by time
$H_{j}$, with $ j \le k-1$, on the event $\{H_{j}<S\}$, at least $j$
red balls have been extracted. To see this, we first focus on $H_{1}$,
and prove that by this time, on the event $\{H_{1}<S\}$, at least one
red ball has been picked. Suppose that this is not true; that is,
suppose that we picked 0 red balls by time $H_{1}$. The reader can
check from our embedding that this implies that
\[
\min_{0 \le t \le S} B_{t}> F- F_{k} - 1/f_{1}
= F_{1} - F_{k}.
\]
This would imply that $H_{1}>S$ contradicting our hypothesis. By
reiterating the same reasoning we get that the statement holds true for
any $j \le k$.

Define
\[
M_{j} \df \{\mbox{after time } H_{j-1}, \mbox{ the process }
\mathbf{B} \mbox{ reaches } F_{j+1} \mbox{ before it hits }
F_{j+1}- F_{k} \}.
\]
On $M_{j}$ the Brownian motion, after time $H_{j-1}$, will hit
$F_{j+1}$ before there is a tie in the urn, because $ F_{j+1}- F_{k}
\ge0$, for $j \le k$.
Next we prove that for any $j \in\{1,2, \ldots, k\}$, if $M_{j}$
holds, then only a finite number of red balls are extracted, that is,
$M_{j} \subset D_{2}$. We split this proof into two parts: we first
prove that $ M_{j} \cap\{ S \le H_{j-1}\} \subset D_{2}$ and then $
M_{j} \cap\{ S > H_{j-1}\} \subset D_{2}$.
In order to prove the first inclusion, recall that under $\mathbb
{Q}^{\ssup k}$ the Brownian motion starts at $F - F_{k}$. This implies
that if $S \le H_{j-1}$, then infinitely many balls will be extracted
before the Brownian motion hits $F_{j} -F_{k}$. As $ F - F_{k} > F_{j}
-F_{k} >0$, we have that infinitely many balls will be extracted before
$\mathbf{B}$ hits $0$, that is, before a tie. This implies that
$B_{S}>0$, which in turn implies $ M_{j} \cap\{ S \le H_{j-1}\}
\subset D_{2}$.

Next we prove that $ M_{j} \cap\{ S > H_{j-1}\} \subset D_{2}$. On the
set $ M_{j} \cap\{ S > H_{j-1}\} $, by time $H_{j-1}$ the number of
red balls extracted is at least $j-1$. This implies that
%
\begin{equation}
\label{punto2} B_{S} \le \sum_{j=1}^{\infty}
\bigl(1/f(j)\bigr) - \sum_{t=1}^{j-1}
\bigl(1/f(t)\bigr)= F_{j} \qquad\forall k \ge n.
\end{equation}
This is a consequence of \eqref{punto} and the fact that $n-\phi(n)$ is
a nondecreasing random sequence, and if $n-\phi(n) = j-1$ for some
$n$, then $\lim_{n \ti} n-\phi(n)\ge j-1$. Let
\[
V_{1} \df\inf\{ m_{n} \dvtx m_{n} >
H_{j-1} \mbox{ and } B_{m_{n}} - B_{m_{n-1}}>0\},
\]
that is, the first time after $H_{j-1}$ that a white ball is extracted.
The stopping time $V_{1}$ could be infinite. Next we prove that on
$M_{j}$ the random time $V_{1}$ is a.s. finite. Recall that $H_{j-1}$
is the first time that the process $\mathbf{B}$ hits $ F_{j} -F_{k}$,
and that $ 0 < F_{j} -F_{k} < F - F_{k}$. This implies that by time
$H_{j-1}$ the number of white balls generated by the Brownian motion,
plus the initial $k$, overcomes that of the red ones. On $M_{j}$, after
time $H_{j-1}$, the process will hit $F_{j+1}$ before it hits $0$. This
implies that $V_{1}<\infty$ a.s. on $M_{j}$. In fact if no white balls
are extracted after time $H_{j-1}$ the process would hit $0$ before it
hits $F_{j+1}$ giving a contradiction.
Moreover on $M_{j}$, we have that $B_{V_{1}}>0$, hence by time $V_{1}$
the white balls are still ahead with respect the red ones.
We can repeat the same reasoning with
\[
V_{2} \df\inf\{ m_{n} \dvtx m_{n} >
V_{1} \mbox{ and } B_{m_{n}} - B_{m_{n-1}}>0\},
\]
to argue that $V_{2}$ is a.s. finite and by time $V_{2}$ the white
balls are still in advantage. By reiterating this argument, we get that
only finite many red balls will be extracted, because each $V_{i}$
occurs before a tie, a.s.
Hence $D_{2}$ holds when $M_{j}$ holds.
This implies that $D_{2}^{c} \subset M_{j}^{c}$ for each $j \in\{
1,2,\ldots, k-1\}$. If $\bigcap_{j=1}^{k-1} M_{j}^{c}$ holds, then either
$\{B_{S}>0\}$ holds or $\{H<S\}$ holds.
If the latter event holds, independently of the past, the probability
that only finitely many white balls are picked is exactly 1/2, by symmetry.
Moreover, the events $M_{j}$ are independent, because they are
determined by the behavior of disjoint increments of the Brownian
motion. By the standard ruin problem for this process, we have that
%
\begin{equation}
\label{fintou1} \mathbb{Q}^{\ssup k}(M_{j}) =
\frac1{1+f(j)F_{k}}.
\end{equation}
%
We get
\[
\mathbb{Q}^{\ssup k}(D_{1}) = \mathbb{Q}^{\ssup k}
\bigl(D_{2}^{c}\bigr) \le\frac 12 \prod
_{j=1}^{s} \biggl(1 - \frac1{1+f(j)F_{k}}
\biggr).
\]
\upqed\end{pf*}

\section{\texorpdfstring{Proofs of Theorem \protect\ref{testo1} and Corollary  \protect\ref{rf}}
{Proofs of Theorem 1.7 and Corollary 1.8}}
\label{lead}
\mbox{}
\begin{pf*}{Proof of Theorem~\ref{testo1}}
Notice that \textit{Lead} must be a vertex of $\mathcal{T}_{1}$. Under the
assumptions of the theorem, the probability that $\inf_{i}x^{*}_{i}> M$
is smaller or equal to the probability that $x^{*}_{1}>M$. The latter
probability is bounded as follows:
\[
\mathbb{P}\Biggl(\sum_{j=1}^{\infty}W^{\ssup1}_{j}/f
\bigl(N_{1}(j) \bigr) > M\Biggr) \le\exp\bigl\{ -(1-p)^{2}
a_{1} (M - 3 F)\bigr\}.
\]
We set $C_{1} = \exp\{3(1-p)^{2}F\}$ and $C_{2} = (1-p)^{2} a_{1}$,
where $a_{1} = \inf_{k \ge1} f(k)$.
In virtue of \eqref{cmlf}, the probability that all the vertices at
level $n$ are good is at least
\[
1 - m^{n} \inf_{r >1} \e^{-c_{n}(r,M) n} + r^{-n},
\]
where $c_{n}(r,M)$ were introduced at the end of the proof of Lemma \ref
{genac}, and $m$ was introduced in \eqref{mo}. Recall that $g_{n}$ is
the set of the vertices of $\Gcal$ at level $n$. Moreover, recall that
$G_{n} = \bigcup_{j \ge n} g_{j}$. We have
\begin{eqnarray*} \mathbb{P}(\mathrm{Lead} \in G_{n})
&\le&\mathbb{P} \Bigl(\Bigl\{\inf_{i} x^{*}_{i}>M
\Bigr\}\cup\{\mbox{at least one vertex in $g_{n}$ is not good}\} \Bigr)
\\[-2pt]
&\le& C_{1}\e^{-C_{2} M} + m^{n} \inf_{r >1}
\bigl(\e^{-c_{n}(r,M) n} + r^{-n}\bigr).
 \end{eqnarray*}
\upqed\end{pf*}

\begin{pf*}{Proof of Corollary~\ref{rf}}
Set $i(1) = \tau_{2}$ and define recursively $i(n) = \inf\{j > i(n-1)
\dvtx R^{\ssup1}_{j} =1\}$. Notice that $i(k) \ge k$. If a vertex $\nu
$ of $\Gcal$ belongs to $g_{1}$, then we have that $\tau_{\nu} = i(k)$
for some $k$. We have\vspace*{-1pt}
\[
 \mathbb{P} (\mathrm{Lead}=1 ) \ge1 -
\mathbb{E}\biggl[\sum_{j
\in g_{1}} \1_{\{x^{*}_{j} < x^{*}_{1}\}}\biggr]
- \mathbb{P}(\mathrm{Lead} \in G_{2}).
\]
We bound the last probability in the previous expression using
Theorem~\ref{testo1}.
Order the groups at level one, starting from the smaller. As $i(k) \ge
k$, we have that by the time the $k$th group at level 1 is created,
there are at least $k$ balls in urn 1. Hence, using Theorem \ref
{estur}, we get
\[
\mathbb{E}\biggl[\sum_{j \in g_{1}} \1_{\{x^{*}_{j} < x^{*}_{1}\}}
\biggr] \le \sum_{k=1}^{\infty} \frac12 \prod
_{\ell=1}^{k-1} \frac{f(\ell)F_{k}}{1 +
f(\ell) F_{k}}.\vspace*{-1pt}
\]
\upqed\end{pf*}

\begin{appendix}
\section*{Appendix}\label{app}
Fix two real numbers $r$ and $w$, and two sequences of positive real
numbers $\{W(k), k \ge w\}$ and $\{R(i), i \ge r\}$. Suppose we have
an urn with
$w$ (resp., $r$) white (resp., red) balls. If at step $n\ge0$ there are
exactly $j$ white balls, with $ n-j \ge0 \ge w -j$, then the probability
to pick a white ball is
\[
\frac{W(j)}{W(j) + R(n-j+w)}.
\]
If a white (resp., red) ball is picked, at time $n+1$ the composition of
the urn becomes $j+1$ (resp., j) white balls and $ n-j+w$ (resp.,
$n-j+w+1$) red ones. Denote by
\begin{eqnarray*}
A^{c}_{R} &\df& \{ \mbox{ the number of red balls in the
urn goes to $\infty$ as $n \ti$}\},
\\
A^{c}_{W} &\df& \{ \mbox{ the number of red balls in the
urn goes to $\infty$ as $n \ti$}\}.
\end{eqnarray*}
Let $\mathbb{Q}$ be the measure describing the dynamics of this urn.
We have the following result, due to Herman Rubin; see the Appendix in
\cite{BD1990}.\vadjust{\goodbreak}
%
\begin{theorem}[(H. Rubin)]\label{hrubin} We have the following 3 cases:
\begin{longlist}[(iii)]
\item[(i)] If $\sum_{k=w}^{\infty}  (W(k) )^{-1} =\infty$ and
$\sum_{k=r}^{\infty}  (R(k) )^{-1} =\infty$, then both the number
or red balls and the number of white balls in the urn goes to $\infty$,
a.s., as $ n \ti$.
\item[(ii)] If $\sum_{k=w}^{\infty}  (W(k) )^{-1} <\infty$ and
$\sum_{k=r}^{\infty}  (R(k) )^{-1} =\infty$, then
\[
\mathbb{Q}(A_{R}) =1.
\]
\item[(iii)] If $\sum_{k=w}^{\infty}  (W(k) )^{-1} <\infty$ and
$\sum_{k=r}^{\infty}  (R(k) )^{-1} <\infty$, then
\[
\mathbb{Q}(A_{R}) +\mathbb{Q}(A_{W}) =
\mathbb{Q}(A_{R} \cup A_{W}) =1,
\]
and both $\mathbb{Q}(A_{R})$ and $\mathbb{Q}(A_{W})$ are strictly positive.
\end{longlist}
\end{theorem}
\end{appendix}

\section*{Acknowledgments}
We thank two anonymous referees for
their suggestions. Moreover we thank Peter M\"orters for
helpful discussions, Patrick Lahr for spotting a few typos in an early
version and Roman Koteck\'y for pointing out to us the reference \cite
{Pric}, which was one of the first to introduce preferential attachment
schemes.

\printaddresses

\end{document}